\documentclass{amsart}
\usepackage[T1]{fontenc}
\usepackage{amssymb}
\usepackage[hyperref]{degt}

\def\2{\color{red}}

\def\newcs#1#2{\expandafter\gdef\csname#1\endcsname{#2}}

\def\BB{$B2$}
\def\Nmax#1{N_{#1}}

\def\G{\frak G}
\def\bS{\bold S}
\def\lattice{\bS\sb{h}}
\def\Sperp{\bold{T}}
\def\(#1){(#1)}
\def\(#1){[#1]}
\def\L{\bold L}
\def\th{\tilde h}

\def\Golay{\Cal{G}}
\def\CC{\Cal{C}}
\def\CK{\Cal{K}}
\def\CO{\Cal{O}}
\def\KC{\bar\CO}
\def\KK{\bar\CK}
\def\KC{\CO_*}
\def\KK{\CK_*}
\let\dm\Omega
\let\graph\Gamma

\let\eq\sim
\let\Eq\approx
\def\cl#1{[#1]}
\def\Cl#1{[\![#1]\!]}
\def\sd{\mathbin\vartriangle}

\def\supp{\operatorname{supp}}
\def\Fn{\operatorname{Fn}}
\def\Fano{\Cal{F}}
\def\irr{^*}

\def\EQ#1){}%
\let\NE\to

\def\bX{\bold X}
\def\bY{\bold Y}

\def\BS{\GROUP{BB}}
\def\BS{\operatorname{K}}

\def\OGplus{\OG^{+\!}}
\def\Eplus{E_+}

\def\famdm#1{\tilde{#1}}
\def\famR{\Cal R}
\def\famX{\Cal X}
\def\famY{\Cal Y}
\def\famZ{\Cal Z}
\def\famXdm{\famdm\famX}

\def\famB{\Cal B}
\def\famM{\Cal M_{20}}
\def\br{\famdm r}
\def\bc{\famdm c}
\def\bn{\famdm n}
\def\rc{\operatorname{rc}}

\let\rs\tau

\def\gk{\frak k}
\def\go{\frak o}
\def\gr{\frak r}
\def\gs{\frak s}
\def\gu{\frak u}
\def\ggu{\gu'}
\def\gv{\frak v}
\def\tu{\tilde u}
\def\hu{\hat u}

\def\dd{\operatorname{d}}

\let\hh\hbar
\def\ds#1#2{\mathopen\|#1/#2\mathclose\|}

\def\cn#1-#2{#1\sb{#2}}
\def\ln#1-#2{#1^\star\sb{#2}}
\let\lvertex\bullet

\def\dual{^\vee}

\def\Sym{\operatorname{Sym}}

\def\CN#1{#1}
\def\LN#1{#1^\star}
\def\STRAT#1{\ifcase#1 \or
 \LN1\or
 \LN2\or
 \CN3\or
 \CN4\or
 \CN5\or
 \CN6\or
 \CN7\fi}
\def\astrat#1{\raise6pt\hbox{\hypertarget{strat-#1}{}}\STRAT#1}
\def\rstrat#1{\hyperlink{strat-#1}{\STRAT#1}}
\pdef\strat#1{\ifmmode\rstrat#1\else$\rstrat#1$\fi}

\def\rsame#1#2{\hyperlink{same-#1-#2}{${}^{#2}$}}
\def\same#1{\expandafter\ifx\csname same-\SAME-#1\endcsname\relax
 \smash{\llap{\raise9pt\hbox{\hypertarget{same-\SAME-#1}{}}${}^{#1}$\ }}
 \expandafter\gdef\csname same-\SAME-#1\endcsname{#1}
\else\smash{\llap{\rsame\SAME#1\ }}\fi}
\pdef\symplectic(#1,#2){\expandafter\ifx\csname(#1,#2)\endcsname\relax??\else
 \csname(#1,#2)\endcsname\fi}

\def\mrsign#1{$\ifcase#1\or*\or*\fi$}
\def\maxr#1{\llap{\smash{\raise9pt\hbox{\hypertarget{maxreal-#1}{}}}\mrsign#1\ }}
\pdef\maxreal#1{\hyperlink{maxreal-#1}{\mrsign#1}}

\def\tabdefs{%
 \let\sep\
 \let\++%
 \def\(##1,##2){\symplectic(##1,##2)}%
 \def\^##1{\ifx1##1\else\rlap{$^##1$}\fi}%
 \let\s\strat
}

\title{Conics in Kummer quartics}

\author{Alex Degtyarev}

\address{%
Department of Mathematics\\
Bilkent University\\
06800 Ankara, TURKEY}

\email{
degt@fen.bilkent.edu.tr}

\thanks{%
The author was partially supported by the T\"{U}B\DOTaccent{I}TAK grant 118F413%
}



\begin{document}

\begin{abstract}
We classify the configurations of lines and conics in smooth Kummer quartics,
assuming that all $16$ Kummer divisors map to conics.
We show that the number of conics on such a quartic is at most $800$.
\end{abstract}

\maketitle

\section{Introduction}\label{S.intro}

All algebraic varieties considered in this paper are over~$\C$.

Following Bauer~\cite{Bauer:conics},
define $\Nmax{2n}(d)$ as the maximal number of
smooth rational curves of degree~$d$
that can lie on a smooth degree~$2n$ $K3$-surface
$X\subset\Cp{n+1}$.
The numbers $\Nmax{2n}(1)$ are quite well understood (see
\cite{Segre,rams.schuett,DIS,degt:lines,degt:sextics}),
whereas for $d=2$ the only known sharp bound is $\Nmax6(2)=285$
(see~\cite{degt:conics}). In the most interesting case of spatial quartics,
$2n=4$, there are but three sporadic examples, with $352$ (Barth,
Bauer~\cite{Barth.Bauer:conics}), $432$ (Bauer~\cite{Bauer:conics}), and $800$
conics (see my recent paper~\cite{degt:800}, where I also motivate the
conjecture that $\Nmax4(2)=800$).

The construction of~\cite{Barth.Bauer:conics,Bauer:conics} embeds a Kummer
$K3$-surface to~$\Cp3$ so that each of the $16$ Kummer divisors is mapped to
a smooth conic, and \latin{a posteriori} (X.~Roulleau, private
communication),
the quartic found in~\cite{degt:800} is also of this nature.
Therefore, we define a \emph{Barth--Bauer quartic} as a smooth quartic
$X\subset\Cp3$ containing $16$ pairwise disjoint smooth conics.
All Barth--Bauer quartics constitute a $3$-parameter family~$\famB$,
%(throughout the paper, we refer to the dimension modulo $\PGL(4,\C)$),
and the
main goal of this paper is a detailed analysis, in the spirit
of~\cite[\S6]{degt.Rams:octics}, of the configurations of lines and conics on
quartics $X\subset\famB$.

The principal results are stated in \autoref{s.results}, after the necessary
preparation in \autoref{s.notions}. We substantiate the conjecture that
$\Nmax4(2)=800$, as well as a few other speculations motivated
by~\cite{degt:conics,degt.Rams:octics}, and discover plenty of examples of
smooth spatial quartics with many conics, both irreducible and reducible.

\subsection{A few basic concepts}\label{s.notions}
Throughout this paper, a \emph{line} on a smooth quartic surface $X\subset\Cp3$
is a smooth rational curve $l\subset X$ of projective degree~$1$,
whereas a \emph{conic} is a curve $c\subset X$ of projective degree~$2$ and
arithmetic genus~$0$. Thus, \emph{we do not assume a conic irreducible}:
it may split into a pair of intersecting lines.

We associate with~$X$ its graphs $\Fn_1X$ and $\Fn_2X$ of lines and conics,
respectively: their vertices are, respectively, lines and
%\emph{irreducible}
conics on~$X$, and two vertices $c_1,c_2$ are connected by an edge of
multiplicity $c_1\cdot c_2\in\Z$.
We denote by $\Fn\irr_2X\subset\Fn_2X$ the graph of \emph{irreducible}
conics, and the union
%The union
\[*
\Fn X:=\Fn_1X\cup\Fn\irr_2X,
\]
with the vertices further connected as above and colored according to their
%projective
degree,
is called the full \emph{Fano graph} of~$X$. The reducible conics are
recovered from $\Fn X$ as pairs
$\mathrel\lvertex\joinrel\relbar\joinrel\relbar\joinrel\mathrel\lvertex$
of adjacent $1$-vertices.

Given an abstract bi-colored graph~$\graph$, we are interested in the
\emph{equiconical stratum}
\[*
\famX:=\famX(\graph):=\bigl\{X\in\famB\bigm|\Fn X\cong\graph\bigr\}\subset\famB,
\]
provided that it is nonempty. Clearly, a necessary condition is
$\graph\supset\dm$, where $\dm$ is the discrete graph of sixteen $2$-vertices
(here and below, a graph inclusion always refers to an induced subgraph with the induced
coloring), and we actually classify the so-called \emph{relative forms}
of~$\graph$, \ie, pairs $(\graph,\dm)$ considered up to automorphism
of~$\graph$. When considering the absolute and relative groups
\roster*
\item
$\Aut\graph\supset\Aut(\graph,\dm)$ of abstract graph automorphisms,
\item
$\Aut_hX\supset\Aut_h(X,\dm)$ of projective automorphisms, and
\item
$\Sym_hX\supset\Sym_h(X,\dm)$ of projective symplectic automorphisms,
\endroster
we always assume that $X\in\famX$ is a very general member.

\remark\label{rem.Omega}
According to Nikulin~\cite{Nikulin:Kummer},
%(see also \autoref{s.general}),\mnote{to check the ref when written}
a set~$\dm$
of $16$ pairwise disjoint irreducible conics on a quartic~$X$ inherits
from~$X$ a certain intrinsic structure
(the $4$-Kummer structure in \autoref{s.Kummer}), and \latin{a priori} it is not
obvious that this structure is preserved by the full group $\Aut\graph$.
However,
\autoref{cor.structure} states that, indeed, the $4$-Kummer
can be
recovered solely from the graph $\graph=\Fn X$; in particular, it is
safe to define the equivalence of relative forms using
the full automorphism group
$\Aut\graph$.
\endremark

For each equiconical stratum $\famX(\graph)$, one can consider the covering
$\famXdm(\graph)\to\famX(\graph)$ consisting of pairs $(X,\dm)$, where
$X\in\famX$ and $\dm$ is a distinguished unordered set of $16$
pairwise disjoint irreducible conics on~$X$. There is an obvious splitting
\[*
\famXdm(\graph)=\bigsqcup\famXdm_i,\quad \famXdm_i:=\famXdm(\graph,\dm_i),
\]
the union running over all relative forms $(\graph,\dm_i)$ of~$\graph$.
The complex conjugation $X\mapsto\bar X$ induces a well-defined involution on the
set $\pi_0(\famY)$ of connected components of each (sub-)stratum $\famY:=\famX$,
$\famXdm$, $\famXdm_i$, \etc.;
we denote by
\[
%(r,c)
\rc\famY
:=\bigl(r(\famY),c(\famY)\bigr)
\label{eq.(r,c)(Y)}
\]
the numbers
of, respectively, real components of~$\famY$ and pairs of complex conjugate
ones (\ie, respectively, one- and two-element orbits of the conjugation).

\subsection{Principal results}\label{s.results}
The main result of the paper is a complete classification, up to
equiconical deformation, of all Barth--Bauer quartics $X\in\famB$ and pairs
$(X,\dm)$, where $\dm$ is a distinguished unordered set of $16$ pairwise disjoint
irreducible conics on~$X$.
The findings are presented in Tables~\ref{tab.1}, \ref{tab.2},
and~\ref{tab.3} (see \cite{degt:4Kummer.tables} for the full version). Let
$n$ and~$\bn$ be the numbers of isomorphism classes of Fano graphs and
relative forms, respectively, and let $(r,c)$ and $(\br,\bc)$ be the
corresponding total component counts, see~\eqref{eq.(r,c)(Y)}.
Itemizing by the codimension of the strata in~$\famB$, we have:
\roster*
\item
$\codim=0$: $(n;r,c)=(\bn;\br,\bc)=(1;1,0)$, see \autoref{tab.1};
\item
$\codim=1$: $(n;r,c)=(\bn;\br,\bc)=(7;7,0)$, see \autoref{tab.1};
\item
$\codim=2$: $(n;r,c)=(43;43,0)$, $(\bn;\br,\bc)=(47;49,0)$, see \autoref{tab.2};
\item
$\codim=3$: $(n;r,c)=(211;208,189)$, $(\bn;\br,\bc)=(285;231,682)$, see
%\autoref{tab.3}.
\cite{degt:4Kummer.tables}.
\endroster

Apart from a plethora of examples of smooth quartics with many conics, these
data are of rather technical nature.
%In view of the technical nature of these data,
For this reason, we state also a few
immediate
consequences that are of a more general interest.
Proofs, if any, are found in \autoref{S.proofs}.

\theorem[see \autoref{proof.main}]\label{th.main}
The maximal number of conics on a Barth--Bauer quartic is $800$.
Up to projective transformation, there is a unique
Barth--Bauer quartic with $800$ conics\rom: it is the Mukai quartic
\[*
z_0^4+z_1^4+z_2^4+z_3^4+12z_0z_1z_2z_3=0
\]
admitting a faithful symplectic action of the Mukai group~$M_{20}$,
see~\cite{Mukai}.
\endtheorem

A smooth quartic with $800$ irreducible conics was first discovered
in~\cite{degt:800}; then, X.~Roulleau observed (with a reference
to~\cite{Bonnafe.Sarti}) that this quartic must be given by Mukai's
polynomial above, upon which
%almost immediately
B.~Naskr\k{e}cki found explicit
equations of all $800$ conics.
Together with~\cite{degt:800}, \autoref{th.main} substantiates the conjecture
that $800$ is the sharp upper bound on the number of conics on a smooth
quartic surface.

\theorem[see Tables~\ref{tab.1}, \ref{tab.2}, and~\ref{tab.3}]\label{th.lines}
The maximal number of lines on a Barth--Bauer quartic is~$48$. Up to
projective transformation, there are four Barth--Bauer quartics \rom(two real
and two complex conjugate\rom) with $48$ lines, see \autoref{tab.3}.
\done
\endtheorem

Recall that the maximal number of lines on a smooth quartic is $64$ (see
\cite{Segre,rams.schuett}) and, according to~\cite{DIS}, the only quartic
with $64$ lines is Schur's quartic~\cite{Schur:quartics}
\[*
z_0(z_1^3-z_0^3)=z_2(z_3^3-z_2^3).
\]
Remarkably, this quartic, denoted~$\bX_{64}$ in~\cite{DIS}, is also Kummer,
but its Kummer divisors map to lines rather than conics.

The next statement illustrates the speculation (see~\cite{degt:conics}) that
reducible conics do not affect the upper bounds, as
a quartic with many conics has no lines. (A similar phenomenon is
observed in~\cite{degt.Rams:octics}, where we study lines on octic
$K3$-surfaces that are allowed to have singularities:
the presence of exceptional divisors does reduce the upper bound
on the number of lines.)

\theorem[see Tables~\ref{tab.1}, \ref{tab.2}, and~\ref{tab.3}]\label{th.irreducible}
Let $X$ be a Barth--Bauer quartic. Then\rom:
\roster
\item\label{i.560}
if $\ls|\Fn_2X|>560$,
%if $X$ has more than $560$ conics,
then $X$ is a singular $K3$-surface\rom;
\item\label{i.576}
if $\ls|\Fn_2X|>576$,
%if $X$ has more than $576$ conics,
then $X$ has no lines or reducible conics\rom;
\item\label{i.544}
if $\ls|\Fn\irr_2X|>544$,
then $X$ has no lines or reducible conics.
\done
\endroster
\endtheorem

Following the tradition, we also consider the problem over~$\R$, \ie, try to
estimate the number of \emph{real} conics on a \emph{real} Barth--Bauer quartic.
Recall that a \emph{real structure} on a complex surface~$X$ is an
anti-holomorphic involution $\rs\:X\to X$, and a curve $c\subset X$ is called
\emph{real} (with respect to~$\rs$) if $\rs(c)=c$.
In the next statement, we do \emph{not} assume in advance that the real
structure should preserve (any) collection of $16$ Kummer conics, even though
it is the case for the maximizing family.

\theorem[see \autoref{proof.real}]\label{th.real}
The maximal number of \emph{real} conics on a \emph{real} Barth--Bauer quartic is $560$.
There is a unique $1$-parameter
equiconical family \rom(the row marked with a \maxreal1
in \autoref{tab.2}\rom) of real Barth--Bauer quartics with $560$ real
conics\rom; all conics are irreducible and have real points.
\endtheorem

The closure of the $1$-parameter family given by
\autoref{th.real} contains a single (up to projective transformation)
singular quartic (see \maxreal2 in
\autoref{tab.3})
with the same collection
of $560$ real conics and $24$ pairs of complex conjugate ones.

The last two theorems should extend to all smooth quartics; however,
unlike \autoref{th.main}, the bounds given by Theorems~\ref{th.irreducible}
and~\ref{th.real} would need to be corrected.
Some counterexamples, with many reducible conics,
can be found in~\cite{DIS}.
Thus,
\roster*
%\item
%$\bX_{64}$: $720=576+144$;
%\item
%$\bX_{60}^\mathrm{i}$: $640=500+140$;
%\item
%$\bX_{60}^\mathrm{ii}$: $680=510+170$;
%\item
%$\bX_{56}$: $624=440+184$;
%\item
%$\bY_{56}$: $636=448+188$;
%\item
%$\bQ_{56}$: $656=448+208$;
%\item
%$\bX_{54}$: $520=404+116$;
%\item
%$\bQ_{54}$: $568=408+160$;
%\item
%$\bZ_{52}$: $480=384+96$;
\item
$\bX_{64}$ has $\text{$576$ reducible}+\text{$144$ irreducible}=720$ conics,
\cf. \autoref{th.irreducible}\iref{i.576}.
\endroster
In fact, $576$ is the maximal number of reducible conics in a smooth quartic
and $720$ appears to be the correct bound for
\autoref{th.irreducible}\iref{i.576}. As another example,
\roster*
\item
the quartic $\bY_{56}$ with $56$ (the maximal number, see~\cite{DIS}) real
lines given by
\[*
\let\e\epsilon
3z_0^2z_1z_2+3z_1z_2z_3^2-z_1^3z_2-z_1z_2^3+2\e z_0^3z_3-2\e z_0z_3^3=0,
\quad
%\e=\pm\sqrt2.
\e^2=2.
\]
has $\text{$448$ reducible}+\text{$188$ irreducible}=636$ conics,
all real, \cf. \autoref{th.real}.
\endroster
It appears that
Barth--Bauer quartics are not good candidates for maximizing the number of
real conics; this fact is discussed in \autoref{s.real}.
I would conjecture the existence of a smooth real quartic with more than $720$
irreducible real conics.

%Thus,
%\roster*
%%\item
%%$\bX_{64}$: $720=576+144$;
%%\item
%%$\bX_{60}^\mathrm{i}$: $640=500+140$;
%%\item
%%$\bX_{60}^\mathrm{ii}$: $680=510+170$;
%%\item
%%$\bX_{56}$: $624=440+184$;
%%\item
%%$\bY_{56}$: $636=448+188$;
%%\item
%%$\bQ_{56}$: $656=448+208$;
%%\item
%%$\bX_{54}$: $520=404+116$;
%%\item
%%$\bQ_{54}$: $568=408+160$;
%%\item
%%$\bZ_{52}$: $480=384+96$;
%\item
%$\bX_{64}$ has $\text{$576$ reducible}+\text{$144$ irreducible}=720$ conics,
%\cf. \autoref{th.irreducible}\iref{i.576};
%in fact, $576$ is the maximal number of reducible conics in a smooth quartic
%and $720$ appears to be the correct bound for
%\autoref{th.irreducible}\iref{i.576};
%\item
%the quartic $\bY_{56}$ with $56$ (the maximal number, see~\cite{DIS}) real
%lines given by
%\[*
%\let\e\epsilon
%3z_0^2z_1z_2+3z_1z_2z_3^2-z_1^3z_2-z_1z_2^3+2\e z_0^3z_3-2\e z_0z_3^3=0,
%\quad
%%\e=\pm\sqrt2.
%\e^2=2.
%\]
%has $\text{$448$ reducible}+\text{$188$ irreducible}=636$ conics,
%all real, \cf. \autoref{th.real};
%I would conjecture the existence of a smooth real quartic with more than $720$
%irreducible real conics.
%\endroster
%It appears that
%Barth--Bauer quartics are not good candidates for maximizing the number of
%real conics; this fact is discussed in \autoref{s.real}.
%%the true upper bound on the number of real conics,
%the reason is the fact that, in
%a Kummer quartic, the classes of real lines and conics can span a lattice
%of rank at most~$19$ (see \autoref{cor.real.singular}).

\proposition[see \autoref{proof.NS}]\label{prop.NS}
Let $X\in\famX$ be a very general member of
an equiconical stratum of Barth--Bauer quartics.
Then the lattice $\NS(X)$ is generated
over~$\Z$ by the classes of lines and conics on~$X$.
In particular, $\OG_h(\NS(X))=\Aut(\Fn X)$.
\endproposition

\corollary[see \autoref{proof.NS}]\label{cor.structure}
Let $\graph:=\Fn X$ be the full Fano graph of a Barth--Bauer quartic. Then any
discrete subgraph $\dm\subset\graph$ with $16$ $2$-vertices inherits
from~$\graph$ a canonical $4$-Kummer structure \rom(see
\autoref{s.Kummer}\rom).
\endcorollary

\subsection{Contents of the paper}
In \autoref{S.BB}, we describe the homological properties of the Barth--Bauer
quartics, mainly using Nikulin~\cite{Nikulin:Kummer}.
In \autoref{S.strata}, we discuss the equiconical strata and their connectedness
(following Nikulin~\cite{Nikulin:forms} and
Dolgachev~\cite{Dolgachev:polarized}); then, we describe
the construction by means of extensions
\via\ extra lines and conics.

The results of the computation are presented in
\autoref{S.codim1}--\autoref{S.codim3}, where we treat the strata of
codimension $1$ to~$3$, respectively. \autoref{S.proofs} fills
in a few missing details in the proofs.

\subsection{Common notation}\label{s.notation}
Unless stated otherwise, all lattices considered in the paper are
\emph{even}: $x^2=0\bmod2$ for all $x\in L$.
We use $\oplus$ for \emph{orthogonal} direct sums of lattices, as opposed to
$+$ for direct, but not necessarily orthogonal sums. For a lattice~$L$, we
use the following notation:
\roster*
\item
$\det L$ is the determinant of the Gram matrix of~$L$ in any integral basis;
\item
$L$ is \emph{nondegenerate} (resp.\ \emph{unimodular}) if $\det L\ne0$
(resp.\ $\det L=\pm1$);
\item
$\Gs_\pm(L)$ are the inertia indices of the quadratic space $L\otimes\R$;
%\item
%$L$ is \emph{hyperbolic} if $\Gs_+L=1$;
\item
$\OG(L)$ is the orthogonal group of~$L$;
\item
$\OGplus(L)\subset\OG(L)$ is the subgroup preserving the
\emph{positive sign structure}, \ie,
coherent orientation of maximal positive
definite subspaces of $L\otimes\R$;
\item
$L\dual:=\{x\in L\otimes\Q\,|\,\text{$x\cdot y\in\Z$ for all $y\in L$}\}$
is the \emph{dual group};
\item
$\discr L:=L\dual\!/L$ is the \emph{discriminant form} of~$L$
(see $q_L$ in~\cite{Nikulin:forms});
\item
$\discr_pL:=(\discr L)\otimes\Z_p$ is the $p$-primary part for a prime~$p$;
\item
$\Aut(\discr L)$ \etc.\ are the groups of auto-isometries of discriminant forms;
\item
$\dd_L\:\!\OG(L)\to\Aut(\discr L)$ is the canonical homomorphism;
\item
$nL:=L^{\oplus n}$ for an integer $n\ge1$;
\item
$L(q)$, $q\in\Q$, is the same abelian group with the form
$x\otimes y\mapsto q(x\cdot y)$.
\endroster
We denote by $[a]:=\Z u$, $u^2=a$, a rank~$1$ lattice and
use the inline notation
\[
\text{$[a,b,c]$ stands for the lattice $\Z u+\Z v$,
$u^2=a$, $v^2=c$, $u\cdot v=b$}.
\label{eq.T}
\]
to describe lattices of rank~$2$.
We fix a lattice~$\L$ isomorphic to the second homology
\[
\L:=H_2(X;\Z)\cong2\bE_8\oplus3\bU
\label{eq.L}
\]
of a $K3$-surface; here, $\bE_8$ and $\bU\cong[0,1,0]$
are the unique unimodular even
lattices of signature $(0,8)$ and $(1,1)$, respectively.

Given a lattice $\Z\dm$ with a distinguished basis~$\dm$,
we identify subsets $\gs\subset\dm$
with the vectors
$\gs:=\sum e\in\Z\dm$, the summation running over $e\in\gs$.
We will also use the shorthand notation
\[*
\hh:=\tfrac12h\in\Q h,\qquad
\ds\gs\gr:=\tfrac12(\gs\cap\gr)-\tfrac12(\gs\sminus\gr)\in\Q\dm,
\quad
%\text{where $\gr,\gs\subset\dm$},
\gr,\gs\subset\dm,
\]
where $h$ will commonly stand for the class of the hyperplane section.
We use the deprecated notation $\bar\gs:=\dm\sminus\gs$ for
the complement of a subset $\gs\subset\dm$.

On a $K3$-surface~$X$, we identify $(-2)$-curves $c\subset X$ (irreducible or
reducible) with their classes in
%the N\'{e}ron--Severi lattice
$\NS(X)$.
The \emph{transcendental lattice} $\NS(X)^\perp\subset H_2(X;\Z)\cong\L$ is
denoted by~$T(X)$.
A $K3$-surface~$X$ is \emph{singular} if $\rank\NS(X)=20$.

\subsection{Acknowledgements}
I would like to express my gratitude to
Dino Festi,
Bartosz Naskr\k{e}cki,
S{\l}awomir Rams,
Xavier Roulleau, and
Matthias Sch\"{u}tt for a number of fruitful discussions.

\section{The Barth--Bauer family}\label{S.BB}

In this section, we give a detailed description of the homological properties
of an abstract Kummer surface and then, of a generic Barth--Bauer quartic.

\subsection{Abstract Kummer structures\pdfstr{}{ \rm(see~\cite{Nikulin:Kummer})}}\label{s.Kummer.structure}
Consider the set $\dm:=\{1,\ldots,16\}$, denote $\CO_0:=\{\varnothing\}$ and
$\CO_{16}:=\{\dm\}$, and define
an \emph{abstract Kummer structure}
as a collection $\CO_8$ of $30$ eight-element subsets
$\go\subset\dm$ such that $\CO_*:=\bigcup_n\CO_n$ is closed under the
symmetric difference~$\sd$.
(We extend $\CO_n$ and other similar notations
\via\ $\CO_n:=\varnothing$ unless the index $n\in\Z$ has been
mentioned explicitly.)

An abstract Kummer structure gives rise to the sets
\[*
\CC_*:=\bigl\{\gs\subset\dm\bigm|
 \text{$\ls|\gs\cap\go|=0\bmod2$ for all $\go\in\CO_*$}\bigr\},\quad
\CC_n:=\bigl\{\gs\in\CC_*\bigm|\ls|\gs|=n\bigr\}
\]
and equivalence relation
\[*
\text{$\gr\eq\gs$ iff $\gr\sd\gs\in\KC$
 and equivalence classes $\cl\gs$ and $\cl\gs_n:=\cl\gs\cap\CC_n$, $n\in\Z$}.
\]

According to Nikulin~\cite{Nikulin:Kummer}, an abstract Kummer structure
exists and is unique up to permutation. Consider the setwise stabilizer
$\stab\CO_*\subset\Bbb S(\dm)$ and its subgroup $\G_\Go$ acting identically on
%the set of equivalence classes.
$\CC_*/{\eq}$.
Then, the group $\G_\Go\cong(\Z/2)^4$ (see
$\#\symplectic(16,14)$ in \autoref{tab.groups}) acts simply transitively
on~$\dm$, making it a $\G_\Go$-torsor or, equivalently, a dimension~$4$ affine
space over $\Bbb{F}_2$, and $\CO_8$ is the set of all affine hyperplanes.
In this language, $\G_\Go\subset\stab\CO_*$ are, respectively,
the group of translations and
that of all affine linear transformations of~$\dm$.

For an alternative description, consider the extended binary Golay code (see,
\eg, \cite{Conway.Sloane}) $\Golay_*$ on the set $\{1,\ldots,24\}$ and take
for~$\dm$ a codeword of length~$16$. Then,
\[*
\CO_8=\bigl\{\go\bigm|\go\subset\dm,\go\in\Golay_8\bigr\},\qquad
\CC_*=\bigl\{\go\cap\dm\bigm|\go\in\Golay_*\bigr\}.
\]
(Note that $\CC_n=\varnothing$ if $n$ is odd or $n=2,14$.)
The group $\stab\CO_*$ is the restriction to~$\dm$ of its setwise stabilizer
in the Mathieu group $M_{24}=\stab\Golay_*\subset\SG{24}$.

Consider the orthogonal direct sum $\Z\dm:=\bigoplus\Z e$, $e\in\dm$,
$e^2=-2$.
%To shorten the notation, we will freely identify subsets $\gs\subset\dm$
%with the vectors
%$\gs:=\sum e\in\Z\dm$, the summation running over $e\in\gs$.
An abstract Kummer structure on~$\dm$ defines a finite index extension
$\bS(\CO_*)\supset\Z\dm$, \viz. the one generated over~$\Z$ by the
vectors $e\in\dm$ and $\frac12\go\in\Q\dm$, $\go\in\CO_*$.

\lemma[see Nikulin~\cite{Nikulin:Kummer}]\label{lem.bS}
%According to~\cite{Nikulin:Kummer},
The lattice $\bS:=\bS(\CO_*)$ is the only, up to the action of the orthogonal
group
$\OG(\Z\dm)$, finite index extension of $\Z\dm$ such that
\roster
\item\label{i.roots}
$\bS$ has no vectors of square $(-2)$ other than $\pm e$, $e\in\dm$, and
\item\label{i.primitive}
$\bS$ admits a primitive isometric embedding to the lattice~$\L$
given by~\eqref{eq.L}.
\endroster
The original abstract Kummer structure~$\CO_*$ is recovered from
an overlattice $\bS\supset\Z\dm$
%as above
satisfying the two conditions above
\via\ $\CO_*=\{\go\subset\dm\,|\,\frac12\go\in\bS\}$.
\done
\endlemma

The primitive isometric embedding $\bS\into\L$ is unique up to isomorphism,
and the orthogonal complement is
\[*
\Sperp:=\bS_\L^\perp\cong3\bU(2).
%\label{eq.Sperp}
\]
In particular, one has
\[
\text{$u^2=0\bmod4$,\quad $u\cdot v=0\bmod2$\quad for any $u,v\in\Sperp$}.
\label{eq.Sperp}
\]

The next two statements are immediate consequences of~\cite{Nikulin:Kummer}
and the techniques developed in~\cite{Nikulin:forms} (where the notation
$q_L$ is used for the discriminant $\discr L$).
%For the first one, recall
%that one has $s^2=0\bmod4$ for each $s\in\Sperp$, see~\eqref{eq.Sperp},
%and that $\ls|\gs|\bmod4$ is constant within each equivalence
%class~$\cl{\,\cdot\,}$.

\lemma\label{lem.discr}
Assume an abstract Kummer structure $\CO_*$ and a primitive isometric embedding
$\bS:=\bS(\CO_*)\into\L$ fixed. Then, there are canonical bijective isometries
\[*
\kappa\:\Sperp/2\Sperp\to\discr\Sperp
 \to\discr\bS\to\bS/2\bS\to\CC_*/{\eq},
\]
where the quadratic forms on
$\Sperp/2\Sperp$ and $\bS/2\bS$
are
$s\mapsto\frac14s^2\bmod2\Z$ and $\CC_*/{\eq}$ is endowed with the group
law induced by~$\sd$ and the form
$[\gs]\mapsto\frac12\ls|\gs|\bmod2\Z$.
\done
\endlemma

\lemma\label{lem.kernels}
Under the assumptions of \autoref{lem.discr}, consider
another overlattice $\bS\subset N\subset\L$ primitive in~$\L$ and
denote by $\bS^\perp\!:=N\cap\Sperp$
the orthogonal complement of~$\bS$
in~$N$. Then
$\bS^\perp\!$ is primitive in $\Sperp$ and
$N\subset(\bS\oplus\bS^\perp)\otimes\Q$
is generated over $\Z$ by $\bS\oplus\bS^\perp$ and all vectors
%of the form
\[*
\tfrac12(\gs+s),
\]
where $s\in\bS^\perp$ and $\gs\in\kappa(s\bmod2\bS^\perp)$ is a representative.
In fact, it suffices to fix a basis for
%the $\Bbb{F}_2$-vector space
$\bS^\perp\!/2\bS^\perp$ and take for $s$ one representative of each basis
element.
\done
\endlemma

\subsection{Polarized Kummer surfaces\pdfstr{}{ \rm(see~\cite{Nikulin:Kummer})}}\label{s.Kummer}
Let~$X$ be the Kummer surface of a complex torus~$A$.
The set~$\dm$ of the $16$ Kummer divisors on~$X$ is in a bijection with the
set of order~$2$ points in~$A$, which is a subgroup
$(\Z/2)^4\subset A$; thus, $\dm$ acquires from~$X$
a natural abstract Kummer structure.
Alternatively, the primitive hull of $\Z\dm$ in
$\NS(X)$ must be a lattice~$\bS$ as in \autoref{lem.bS}; by the lemma, $\bS$
gives rise to an abstract Kummer structure on~$\dm$.

Now, assume that $X$ is embedded, as a smooth surface, to the projective
space $\Cp{d+1}$ so that each Kummer divisor $e\in\dm$ is mapped to an
\emph{irreducible} conic. Denote by $h\in\NS(X)$ the class of a hyperplane
section and assume that it is not divisible by~$2$ (which is always the case
if $h^2\ne0\bmod8$ or $h^2<16$, see~\cite{Saint-Donat}).
The lattice
\[*
\Z\dm+\Z h\subset\NS(Z),\qquad h^2=2d,\quad
\text{$h\cdot e=2$ for all $e\in\dm$},
\]
%and its primitive hull $N$ in $\NS(X)$.
%The lattice $\Z\dm+\Z h$
splits as
\[
\Z\dm\oplus\Z\th,\qquad\th:=h+\dm,\quad\th^2=h^2+32,
\label{eq.splitting}
\]
and from~\eqref{eq.Sperp} we conclude that $h^2=0\bmod4$.
Then, by \autoref{lem.kernels}, the primitive hull~$N$ of $\Z\dm+\Z h$
in~$\NS(X)$
contains the overlattice
\[
\text{$\lattice:=\lattice(\CK_*)\supset\Z\dm+\Z h$
 spanned by $\bS(\CO_*)$, $h$, and $\tfrac12(h+\gk)$, $\gk\in\CK_*$},
\label{eq.Sh}
\]
where $\CK_*\subset\CC_*$ is a certain equivalence class, $\CK_*\ne\CO_*$,
such that
$2\ls|\gk|=h^2\bmod8$ for $\gk\in\CK_*$.
A choice of this extra class, which is determined by the smooth embedding
$X\into\Cp{d+1}$, is called an \emph{$h^2$-Kummer structure} on~$\dm$.
This choice defines a coarser equivalence relation
\[*
\text{$\gr\Eq\gs$ iff $\gr\sd\gs\in\KC\cup\KK$
 and equivalence classes $\Cl\gs$ and $\Cl\gs_n:=\Cl\gs\cap\CC_n$, $n\in\Z$}.
\]

If $h^2=4$, a \emph{$4$-Kummer structure}
$\CK_*=\CK_6\cup\CK_{10}$ has $16$ sets of size~$6$ and
$16$ sets of size~$10$; it is unique up to $\stab\CO_*$. The discriminant of
the respective lattice $\lattice$ has $2$- and $3$-torsion:
\[
\discr_2\lattice\cong
 \bmatrix0&\frac12\\\frac12&0\endbmatrix\oplus\bmatrix0&\frac12\\\frac12&0\endbmatrix
  \oplus\bigl[\tfrac14\bigr],\qquad
\discr_3\lattice\cong\bigl[\tfrac49\bigr].
\label{eq.discr}
\]
Note that $\CO_*=\{\gk_1\sd\gk_2\,|\,\gk_1,\gk_2\in\CK_*\}$
is determined by $\CK_*$, provided that the set obtained is an abstract Kummer
structure (which is a restriction on~$\CK_*$).

\subsection{Barth--Bauer quartics}\label{s.quartics}
As explained in the previous section, the N\'{e}ron--Severi lattice of a
Barth--Bauer quartic $(X,\dm)$ with a distinguished set~$\dm$ of Kummer
conics must contain the overlattice $\lattice$ as in~\eqref{eq.Sh} defined by
a certain $4$-Kummer structure on~$\dm$. To complete the description of the
generic N\'{e}ron--Severi lattice, we need a realizability criterion.

\definition[\cf. Saint-Donat~\cite{Saint-Donat}]\label{def.admissible}
A hyperbolic overlattice
%$N\supset\lattice\supset\Z\dm+\Z h$ as in~\eqref{eq.Sh}
$N\supset\lattice$, see~\eqref{eq.Sh},
is called
\emph{admissible} if there is no vector $r\in N$ such that either
\roster
\item\label{i.exceptional}
$r^2=-2$ and $r\cdot h=0$ (\emph{exceptional divisor}), or
\item\label{i.isotropic}
$r^2=0$ and $r\cdot h=\pm2$ (\emph{$2$-isotropic vector}), or
\item\label{i.missing}
$r^2=-2$, $r\cdot h=1$, and $r\cdot e<0$ for some $e\in\dm$
(\emph{missing conic}).
\endroster
(The last condition is due to the fact that we insist that each conic
$e\in\dm$ should be irreducible.)
An admissible lattice~$N$ is called \emph{geometric} if it admits a primitive
isometric embedding $N\into\L$, see~\eqref{eq.L}.
\enddefinition

The next statement is well known. It follows from the surjectivity of the
period map for $K3$-surfaces (Kulikov~\cite{Kulikov:periods}) and the results
of Saint-Donat~\cite{Saint-Donat}; \cf. also~\cite{degt:conics} or
\cite{degt.Rams:octics} for an accurate restatement of~\cite{Saint-Donat}
in the homological language.

\proposition\label{prop.existence}
Consider an overlattice $N\supset\lattice\supset\Z\dm+\Z h$, $h^2=4$,
as in~\eqref{eq.Sh}.
Then, there exists a Barth--Bauer quartic $(X,\dm)\subset\Cp3$ such that
\[*
\bigl(\NS(X);\dm,h\bigr)\cong(N;\dm,h)
\]
if and only if $N$ is geometric.
\done
\endproposition

\corollary\label{cor.BB.NS}
Given a Barth--Bauer quartic $(X,\dm)\subset\Cp3$, the primitive hull~$N$
of $\Z\dm+\Z h$ in $\NS(X)$ is the lattice $\lattice$ given by~\eqref{eq.Sh}.
\endcorollary

\proof
According to \autoref{lem.kernels}, the only other option is that the
sublattice $\Z\th$, $\th:=h+\dm$, is not primitive in~$\Sperp$, \ie,
$\th$ is divisible by~$3$. (Recall that $9\mathrel|36=\th^2$.)
But then, for any $\gk\in\CK_{10}$, the vector
$\frac13\th-\frac12(h+\gk)\in N$ would be an
exceptional divisor, see \autoref{def.admissible}\iref{i.exceptional}.
\endproof

Using the description of the nef cone of a $K3$-surface (see, \eg,
\cite[\S8.1]{Huybrechts}), one can easily compute the sets of lines and conics on a
smooth
quartic~$X$
%solely
in terms of the polarized lattice $N:=\NS(X)\ni h$
(\cf.~\cite{degt:conics}):
\[
\aligned
\Fn_n(N,h)&=\bigl\{u\in N\bigm|\text{$u^2=-2$, $u\cdot h=n$}\bigr\},\quad n=1,2,\\
\Fn\irr_2(N,h)&=\bigl\{u\in\Fn_2(N,h)\bigm|
 \text{$u\cdot v\ge0$ for all $v\in\Fn_1(N,h)$}\bigr\}.
\endaligned
\label{eq.Fn}
\]
For a very general Barth--Bauer quartic $X\in\famB$,
using the lattice $\NS(X)=\lattice$
given by \autoref{cor.BB.NS}, this computation shows that $X$ has no lines
(as $h\cdot u=0\bmod2$ for all $u\in\lattice$)
and it has $352$ conics, all irreducible, \viz.
%(see~\cite{Barth.Bauer:conics}), \viz.
\roster*
\item
the $16$ original Kummer conics $e\in\dm$,
\item
$16$ \emph{dual Kummer conics} $e^*:=h-e$, $e\in\dm$, and
\item
$320$ \emph{Barth--Bauer}, or \emph{\BB-conics}
\[
\hh+\ds\gk\gs,\qquad \gs\subset\gk\in\CK_6,\ \ls|\gs|=3.
\label{eq.10-3}
\]
\endroster
This computation agrees with~\cite{Barth.Bauer:conics} and leads us to the
following corollary.

\corollary\label{cor.Sh.conics}
The N\'{e}ron--Severi lattice $\NS(X)=\lattice$ of a very general Barth--Bauer
quartic $(X,\dm)$ is generated over~$\Z$ by the conics on~$X$.
Hence, one has
\[*
\OG_h\bigl(\NS(X)\bigr)=\OG_h(\lattice)=\Aut(\Fn X)
\]
and the $4$-Kummer structure on $\dm\subset\Fn X$ is recovered
from $\Fn X$.
\done
\endcorollary

\subsection{Symmetries}\label{s.symmetries}
It is immediate that, for a very general Barth--Bauer quartic $(X,\dm)$  as in
\autoref{cor.Sh.conics}, one has
\[*
\Aut(\Fn X,\dm)=\G:=\stab(\CK_*\cup\CO_*),\qquad
\Aut(\Fn X)=\G\times(\Z/2),
\]
where $\G\subset\Bbb{S}(\dm)$
is the setwise stabilizer of the $4$-Kummer structure and
the second factor $\Z/2$ is generated by the duality involution
$c\mapsto c^*:=h-c$. (Geometrically, $c\cup c^*$ is the quartic curve cut
by~$X$ on the plane spanned by~$c$.)

\table
\caption{$\G$-orbits on $\CC_n$}\label{tab.orbits}
\def\0{\phantom0}
\centerline{\small\vbox{\halign{&\strut\quad\hss$#$\hss\quad\cr
\noalign{\hrule\kern3pt}
n&\text{even}&\text{odd}&\KC\cup\KK\cr
\noalign{\kern1pt\hrule\kern3pt}
0&&&1\times1\0\cr
4&15\times4\0&20\times4\0\cr
6&15\times16&12\times16&1\times16\cr
8&15\times24&20\times24&1\times30\cr
10&15\times16&12\times16&1\times16\cr
12&15\times4\0&20\times4\0\cr
16&&&1\times1\0\cr
\noalign{\kern1pt\hrule}
\crcr}}}
\endtable
Given a $4$-Kummer structure~$\CK_*$, we subdivide the sets $\gs\in\CC_*$
into \emph{even} and \emph{odd}, according to the parity
$\ls|\gs\cap\gk|\bmod2$ of the intersection with some (any)
$\gk\in\CK_*$. The $\G$-action on~$\CC_*$ preserves parity and both
$\eq$- and $\Eq$-equivalence.
Each nonempty set
$\CC_n\sminus(\KC\cup\KK)$, $n\in\Z$,
consists of exactly two $\G$-orbits: one even and one odd,
see \autoref{tab.orbits} where the counts are given in the form
\[*
(\text{number of $\eq$-classes})\times(\text{size of a class}).
\]

\table
\caption{Symplectic groups $G_\Go$}\label{tab.groups}
\def\0{\phantom0}%
\makeatletter
\centerline{\small\offinterlineskip\vbox{%
\def\\#1(#2,#3){\omit\raise6pt\hbox{\hypertarget{sym(#2,#3)}{}}\sep\hss$#1$\hss\sep
 \protected@write\@auxout{}{\string\newcs{(#2,#3)}{\string\hyperlink{sym(#2,#3)}{#1}}}%
% \expandafter\gdef\csname(#2,#3)\endcsname{\hyperlink{sym(#2,#3)}{#1}}%
 &#2&#3}%
\let\sep\quad
\halign{&\strut\sep\hss$#$\hss\sep\cr
\noalign{\hrule\kern3pt}
\#&\ls|G_\Go|&\text{index}&G_\Go\cr
\noalign{\kern1pt\hrule\kern3pt}
\\21(16,14)&C_4^2\cr
\\39(32,27)&2^4C_2\cr
\\49(48,50)&2^4C_3\cr
\\56(64,138)&\graph_{25}a_1\cr
\\65(96,227)&2^4D_6\cr
\\75(192,1023)&4^2\frak{A}_4\cr
\\81(960,11357)&M_{20}\cr
\noalign{\kern1pt\hrule}
\crcr}}}
\endtable
In view of \autoref{lem.discr} and the global Torelli theorem for
$K3$-surfaces~\cite{Pjatecki-Shapiro.Shafarevich},
the subgroup $\G_\Go\subset\G$ acting identically on $\CC_*/{\eq}$ (see
\autoref{s.Kummer.structure} and $\#\symplectic(16,14)$ in
\autoref{tab.groups}) is the group of symplectic projective automorphisms of
a very general Barth--Bauer quartic~$X$; hence, it acts symplectically on any
Barth--Bauer quartics. For further references,
all groups appearing in this way (for quartics generic
in their respective equiconical strata) are listed in \autoref{tab.groups},
where $\#$ is a reference to the list found in Xiao~\cite{Xiao:Galois}, ``index'' is the
index in the \GAP's small group library, so that
\[*
G_\Go=\mathtt{SmallGroup}(\ls|G_\Go|,\text{\it index}),
\]
and the last column
is a description of the group in the notation of~\cite{Xiao:Galois}.

\section{Equiconical strata}\label{S.strata}

Below, we discuss the connected components of the equiconical strata (see
\autoref{s.components}) and their construction by means of extensions of the
lattice~$\lattice$ (see \autoref{s.extensions}).
In \autoref{s.remarks}, we
outline the algorithms used in the computation.

\subsection{Connected components\pdfstr{}{ \rm(see~\cite{Dolgachev:polarized,Nikulin:forms})}}\label{s.components}
Fix a bi-colored graph $\graph\supset\dm$ and let
$\Fano(\graph):=(\Z\graph +\Z h)/\ker$, where $\Z\graph$ is freely generated
by the vertices $v\in\graph$ and
%the products are
\[*
h^2=4,\quad
v^2=-2,\quad
v\cdot h=\operatorname{color}(v),\quad
u\cdot v=\#\operatorname{edges}(u,v),\quad
u,v\in\graph.
\]
According to \autoref{prop.existence} and~\eqref{eq.Fn}, the original graph~$\graph$ is the Fano
graph of a Barth--Bauer quartic if and only if there is an extension
\[
%N\supset\Fano(\graph)\quad
%\text{such that $N$ is geometric and $\Fn(N,h)=\graph$}.
N\supset\Fano(\graph)\:\quad
\text{$N$ is geometric, $[N:\Fano(\graph)]<\infty$, and $\Fn(N,h)=\graph$}.
\label{eq.N-extension}
\]
Furthermore,
\[*
\famX(\graph)=\bigsqcup\famX(N),
\]
the union running over
the set of all
such extensions, regarded up to abstract lattice isomorphism preserving~$h$.

Thus, we assume that there is a primitive isometry $\iota\:N\into\L$. Given a
subgroup $G\subset\OGplus(N)$,
%two isometries $\iota_i\:N\into\L$, $i=1,2$, are
%said to be \emph{$G$-equivalent} if there exists a pair
%$\Ga\in G$, $\Gb\in\OGplus(\L)$ such that $\iota_2\circ\Ga=\Gb\circ\iota_1$.
%Similarly, we define
%$G$-equivalence of (anti-)isometries
%$\psi_i\:\discr N\to\discr T$ by requiring that there should exist a pair
%$\Ga\in G$, $\Gb\in\OGplus(T)$ such that
%$\psi_2\circ\dd_N(\Ga)=\dd_T(\Gb)\circ\psi_1$.
introduce the following definitions:
\roster*
\item
two primitive isometries $\iota_i\:N\into\L$, $i=1,2$, are
\emph{$G$-equivalent} if there exists a pair
$\Ga\in G$, $\Gb\in\OGplus(\L)$ such that $\iota_2\circ\Ga=\Gb\circ\iota_1$;
\item
two (anti-)isometries
$\psi_i\:\discr N\to\discr T$ are $G$-equivalent if there exists a pair
$\Ga\in G$, $\Gb\in\OGplus(T)$ such that
$\psi_2\circ\dd_N(\Ga)=\dd_T(\Gb)\circ\psi_1$.
\endroster
Note the usage of $\OGplus$ rather than the more conventional~$\OG{}$:
%it is more suitable for our purposes.
this approach excludes anti-isomorphisms of $K3$-surfaces, \cf.
\autoref{prop.absolute}\iref{abs.real} below.

\proposition[\cf. Nikulin~\cite{Nikulin:forms}]\label{prop.Nikulin}
A primitive isometry $\iota\:N\to\L$ gives rise to a bijective anti-isometry
$\psi:\discr N\to\discr N^\perp$\rom; as a consequence, the genus
of
the lattice
$N^\perp$ is independent of~$\iota$, depending on~$N$ only.

The $G$-equivalence classes of isometries $\iota\:N\into\L$ are in a canonical
bijection with the following sets of data\rom:
\roster
\item\label{Nikulin.T}
a lattice \rom(isomorphism class\rom)
$T\in\operatorname{genus}(N^\perp)$, and
\item\label{Nikulin.psi}
a $G$-equivalence class of anti-isometries
$\psi\:\discr N\to\discr T$.
\done
\endroster
\endproposition

For the isometry $\iota\:N\into\L$ corresponding to a pair $(T,\psi)$ as in
\autoref{prop.Nikulin}, we have
\[
\bigl\{g\in\OG(\L)\bigm|g(N)=N,\ g|_N\in G_N,\ g|_T\in G_T\rom\}
 =G_N\times_{\Aut(\discr T)}G_T,
\label{eq.aut.group}
\]
where $G_N\subset\OG(N)$, $G_T\subset\OG(T)$ is a given pair of subgroups
%of $\OG(\L)$ consisting of the
%automorphisms that preserve~$N$ and restrict to a given pair of subgroups
%$G_N\subset\OG(N)$ and $G_T\subset\OG(T)$ is
and the amalgamated product is over the diagram
\[*
G_N\to\OG(N)\overset{\psi\circ{\dd_N}}\longto\Aut(\discr T)
 \overset{{\dd_T}}\longleftarrow\OG(T)\leftarrow G_T.
\]
(In more formal notation,
$G_N$ and $G_T$ should first be replaced with the respective pull-backs of
the intersection
of their images in $\Aut(\discr T)$.)

%Following~\cite{Dolgachev:polarized},
%We define
%Define (\cf.~\cite{Dolgachev:polarized})
An \emph{$N$-polarized} (more precisely, $(N,h)$-polarized)
$K3$-surface
is defined
as a pair $(X,\Gf)$, where $X$ is a $K3$-surface and
$\Gf\:N\into\NS(X)$ is a primitive isometry such that $\Gf(h)$ is ample.
Two $N$-polarized $K3$-surfaces $(X_i,\Gf_i)$, $i=1,2$, are \emph{isomorphic}
if there exists an isomorphism $f\:X_1\to X_2$ such that
$\Gf_2=f_*\circ\Gf_1$. An $N$-polarized $K3$-surface $(X,\Gf)$
defines a $\{1\}$-equivalence class $\Gf\:N\into\NS(X)\into\L$ of primitive isometries
 and, according to
Dolgachev~\cite{Dolgachev:polarized},
assuming that $N$ satisfies condition~\iref{i.exceptional} in
\autoref{def.admissible} (so that $h$ has a chance to be ample),
each such class gives rise to a coarse
moduli space
\[*
\famZ(N\into\L),
\]
which is a connected quasi-projective variety of
dimension $20-\rank N$.
If $N$ is also admissible, the subspace
\[*
\famZ^\circ(N\into\L):=\bigl\{X\in\famZ(N\into\L)\bigm|
 \text{$X\subset\Cp3$ is smooth, $\Fn(X)=\Fn(N,h)$}\bigr\}
\]
is the complement of the countable collection of divisors
\[*
\bigl\{X\in\famZ(N\into\L)\bigm|\NS(X)\ni r\bigr\},\qquad
r\in\L\sminus N,
%\ r^2=-2;
\]
where $r^2=-2$ and $\ls|r\cdot h|\le2$ or $r^2=0$ and $r\cdot h=2$,
\cf.~\eqref{eq.Fn} and \autoref{def.admissible};
hence, it is still connected. Finally, getting rid of the $N$-marking~$\Gf$,
we have
\[*
\famX(N)/\!\PGL(\C,4)=\left(\bigsqcup\famZ(N\into\L)\right)\!\big/\!\OG_h(N),
\]
where the union
%in the right hand side
runs over all $\{1\}$-isomorphism
classes of isometries.
Since $\PGL(\C,4)$ is connected, we arrive at the following well known
statement.

\proposition\label{prop.absolute}
Fix an extension $N\supset\Fano(\graph)$ as in~\eqref{eq.N-extension},
regarded up to abstract lattice isomorphism preserving~$h$,
and let $G:=\OG_h(N)$. Then\rom:
\roster
\item\label{abs.dim}
the stratum $\famX(N)\subset\famB$ is of pure codimension $\rank N-17$\rom;
\item\label{abs.components}
the connected components of $\famX(N)$ are in a natural bijection with the
$G$-equivalence classes of primitive isometries $N\into\L$, \cf.
\autoref{prop.Nikulin}\rom;
\item\label{abs.real}
a component is real if and only if, in \autoref{prop.Nikulin}, there is an
isometry $\Gb\in\OG(T)\sminus\OGplus(T)$ such that
$\dd_T(\Gb)\in\psi\circ\dd_N(G)$.
\done
\endroster
\endproposition

A similar argument shows that
\[*
\famXdm(\graph,\dm)=\bigsqcup\famXdm(N,\dm),
\]
where the union runs over all finite index extensions $N\supset\Fano(\graph)$ as
in~\eqref{eq.N-extension} regarded up to lattice isomorphism preserving~$h$
and~$\dm$ as a set, and then the connected components of each $\famXdm(N,\dm)$ are
described by an analogue of \autoref{prop.absolute}, with $N$ replaced
with a pair $(N,\dm)$
%,
%regarded up to lattice isometry preserving~$h$ and $\dm$ as a set,
and, respectively, $G$ changed to $\OG_h(N,\dm)$.

\subsection{Extensions of $\lattice$}\label{s.extensions}
Paraphrasing~\eqref{eq.N-extension}, strata of codimension~$r$
are described by geometric extension,
necessarily primitive, $N\supset\lattice$ of corank~$r$ that are
%that may serve as the
%N\'{e}ron--Severi lattices of Barth--Bauer quartics~$X$. Taking for~$X$ a very
%general member of the equiconical stratum
%$\famX\bigl(\Fn(N,h)\bigr)$, we can assume that $N$ is
rationally
generated by lines and conics, \ie, finite index extensions of
lattices of the form
\[
\lattice\(u_1,\ldots,u_r):=\lattice+\Z u_1+\ldots+\Z u_r,
\label{eq.Sh()}
\]
where $u_1,\ldots,u_r$ are extra lines or conics linearly
independent over~$\lattice$.
To complete the description of~\eqref{eq.Sh()}, we need to specify,
for $i,j=1,\ldots,r$,
\roster
\item\label{ex.u.h}
the intersections $u_i\cdot h=1$ or~$2$ if $u_i$ is a line or conic,
respectively,
\item\label{ex.u.e}
the intersections $u_i\cdot e$ for all $e\in\dm$, and
\item\label{ex.Gram}
the Gram matrix $[u_i\cdot u_j]$, where $u_i^2=-2$ for all~$i$.
\endroster
Geometrically,
it is obvious that, given two lines~$l_1$, $l_2$, a conic~$c$, and an
\emph{irreducible} conic~$e$, one has
\[*
l_1\cdot l_2\in\{0,1\},\qquad
l_1\cdot e\in\{0,1,2\},\qquad c\cdot e\in\{0,1,2,4\}
%\label{eq.intr}
\]
and $c\cdot e=4$ if and only if $c+e=h$, in which case
$c$ and~$e$ are dependent
over~$\lattice$. It follows that
%, in item~\eqref{ex.Gram} above, one has
%$u_i\cdot u_j\in\{0,1,2\}$ for $i\ne j$ and
item~\iref{ex.u.e} in the above list is
determined by the \emph{supports}
\[
\supp_nu_i:=\bigl\{e\in\dm\bigm|e\cdot u_i=n\bigr\}\subset\dm,\quad n=1,2,
\label{eq.support}
\]
which, for each~$u_i$, are two disjoint subsets of~$\dm$.
Denoting, for short, $\gu_i:=\supp_1u_i$ and $\ggu_i:=\supp_2u_i$, we have
\[
%\text{$\gu\in\CC$},\qquad
 \text{$\gu_i\in\CC_*$ is even if $u_i$ is a conic},\qquad
 \text{$\gu_i\in\CC_*$ is odd if $u_i$ is a line}
\label{eq.parity}
\]
(to make sure that $u_i\cdot v\in\Z$ for each $v\in\lattice$) and,
if $u_i$ is a conic,
the dual conic $u_i^*:=h-u_i$ has the supports
\[
\supp_1u_i^*=\gu_i,\qquad
\supp_2u_i^*=\smash{\overline{\gu_i\cup\ggu_i}}.
%\dm\sminus(\gu_i\cup\ggu_i).
\label{eq.duality}
\]

\convention\label{conv.lines}
We will construct the lattice $\lattice\(u_1,\ldots)$ inductively, adding one
extra line/conic at a time.
To reduce the overcounting, we will always assume that
\roster*
\item
all extra lines are added first, followed by extra conics, and
\item
the number of lines in the generating set is maximal possible.
\endroster
\endconvention

\lemma\label{lem.convention}
Under \autoref{conv.lines}, the following further restrictions hold\rom:
\roster
\item\label{conv.irreducible}
all extra conics are irreducible, and
\item\label{conv.Gram}
for each pair $u_i\ne u_j$ of distinct extra generators, one has
\[*
\text{$u_i\cdot u_j\in\{0,1,2\}$ if both $u_i$, $u_j$ are conics},\qquad
\text{$u_i\cdot u_j\in\{0,1\}$ otherwise}.
\]
\endroster
\endlemma

\proof
If a conic $u=l_1+l_2$ is reducible in (a finite index extension of) $N+\Z u$,
then
at least one of the lines~$l_1$, $l_2$ is not in~$N$ and can replace~$u$ in the
generating set. For the second statement, in view of the first one, the only
value that needs to be ruled out is $l\cdot c=2$ for a line~$l$ and
conic~$c$. In this case, $l':=h-c-l$ is also a line and $l$, $l'$ generate
the same extension as $l$, $c$.
\endproof

\subsection{Technical remarks}\label{s.remarks}
We conclude this section with a discussion of a few tricks used in
the computation.

On a case-by-case basis, we observe that any extension $N$ as
in~\eqref{eq.N-extension} is necessarily trivial, \ie,
$N=\Fano\bigl(\Fn(N,h)\bigr)$
(\cf. \autoref{prop.NS}). More precisely, when analyzing a corank~$1$
extension
\[*
S_r:=\lattice\(\ldots,u_r)\supset S_{r-1}:=\lattice\(\ldots),
\]
we observe that any proper finite index extension $N\supset S_r$ either fails
to be admissible or has a line or conic $u'_r\notin S_r$ (sometimes violating
\autoref{conv.lines}), so that we can replace $S_r$ with
$\lattice\(\ldots,u'_r)$. Note also that, by induction, we can assume
$S_{r-1}$ primitive in~$N$; therefore,
it suffices to consider further extensions generated by
isotropic vectors of the form
\[*
(\Ga u_r+v)\bmod S_r,\qquad \Ga\in\Q\sminus\Z,\ v\in S_{r-1};
\]
this reduces the amount of computation.

This fact has multiple consequences. First, in \autoref{prop.absolute} (and
its analogue for relative forms) we can take for~$N$ the original lattice
$\lattice\(\ldots,u_r)$. Second, there are natural isomorphisms
\[*
\OG_h(N)=\Aut\Fn(N,h),\qquad
\OG_h(N,\dm)=\Aut\bigl(\Fn(N,h),\dm\bigr),
\]
so that these orthogonal groups used in \autoref{prop.Nikulin} can easily be
computed by the \texttt{Digraphs} package in \GAP~\cite{GAP4}. Finally, the
unevoidable overcounting is easily eliminated by identifying the strata with
isomorphic graphs~$\graph$ or pairs $(\graph,\dm)$.

When computing the number of connected components, \autoref{prop.Nikulin} is
used literally if $\rank N=20$, \ie, the stratum in question has
codimension~$3$ and consists of
singular $K3$-surfaces: in this case, $T$ is positive definite of
rank~$2$ and the genera of such forms and their orthogonal groups are easily
computed, see, \eg, Gauss~\cite{Gauss:Disquisitiones}.
In all other cases, with~$T$ indefinite and
$\OG(T)$ typically infinite, we use the results
of Miranda--Morrison~\cite{Miranda.Morrison:book},
which combine both steps~\iref{Nikulin.T}
and~\iref{Nikulin.psi} in \autoref{prop.Nikulin} into a single homomorphism
\[*
d^\perp\:\!\OG(N)\to\Eplus(T)
\]
%where $\Eplus(T)$ is a certain elementary abelian $2$-group
to a certain elementary abelian $2$-group $\Eplus(T)$
%easily computable
computed
in terms of
%the genus of~$T$.
$\operatorname{genus}(T)$,
so that
%Then,
the $G$-equivalence classes of primitive
isometries $N\into\L$ are in a canonical bijection with $\Eplus(T)/\Im d^\perp$.
We omit the details, which are purely algorithmic.

%Explain the count and move \autoref{ex.forms} here.\mnote{\todo:}

%Furthermore, if $l\cdot c'=2$, then $l':=h-c'-l$ is another line,
%and if $c\cdot c'=4$, then $c=h-c'$; in both extreme cases, $l$ and $l'$
%(respectively, $c$) lie in the plane of~$c'$.

\section{Strata of codimension $1$}\label{S.codim1}

\table
\caption{Strata of codimension $\le1$ \rm(see~\autoref{s.list.1})}\label{tab.1}
\centerline{\vbox{\small\offinterlineskip
\tabdefs
\def\frac#1#2{#1/#2}%
\tabskip0pt plus 1fil
\halign{\strut\sep\hss$#$\hss\sep\tabskip0pt&
 \strut\sep\hss$#$\hss\sep&
 \strut\sep\hss$#$\hss\sep&
 \strut\sep\hss$#$\hss\sep&
 \strut\sep\hss$#$\hss\sep&
 \strut\sep\hss$#$\hss\sep&
 \strut\sep\hss$#$\hss\sep&
 \strut\sep\hss$#$\hss\sep&
 \strut\sep\hss$#$\hss\sep&
 \strut\sep\hss$#$\hss\sep&
 \strut\sep\hss$#$\hss\quad\tabskip0pt plus 1fil\cr
\noalign{\hrule\kern3pt}
\text{Name}&\text{Patterns}&\Gd_2^2&\Gd_3&\text{Lines}&\text{Conics}&\ls|G|&
 i_\dm&G_\Go&\ls|{\det}|&(r,c)\cr
\noalign{\kern1pt\hrule\kern3pt}
\text{open}&&&&&352&23040&2&\(16,14)\^1&576&(1,0)\EQ(1,0)\cr
\noalign{\kern1pt\hrule\kern3pt}
\astrat1&\ln4-0, \ln12-0&\frac14&\pm1& 8& 16\+336& 1152& 2& \(16,14)\^1& 400& (1,0)\EQ(1,0)\cr
\astrat2&\ln6-0, \ln10-0&\frac54&\pm4& 32& 160\+192& 1920& 2& \(16,14)\^1& 208& (1,0)\EQ(1,0)\cr
\noalign{\kern1pt\hrule\kern3pt}
\astrat3&\cn4-0, \cn10-0&1&\pm3& & 392& 3072& 4& \(16,14)\^1& 576& (1,0)\EQ(1,0)\cr
\astrat4&\cn6-0, \cn8-0&1&\pm2& & 432& 3072& 4& \(16,14)\^2& 448& (1,0)\EQ(1,0)\cr
\astrat5&\cn12-0&1&\pm1& & 360& 1536& 2& \(16,14)\^1& 832&(1,0)\EQ(1,0)\cr
\astrat6&\cn12-1&0&\pm4& & 384& 3072& 4& \(16,14)\^1& 640&(1,0)\EQ(1,0)\cr
\astrat7&\cn12-2&1&0& & 400& 3072& 2& \(16,14)\^2& 576& (1,0)\EQ(1,0)\cr
\noalign{\kern1pt\hrule}
\crcr}}}
\endtable

The goal of this section is compiling \autoref{tab.1} listing
the seven equiconical strata of codimension~$1$.
(The first row shows the open stratum of
codimension~$0$.)
%The notation used
%in the table and the meaning of
%its columns are explained in \autoref{s.list.1} below.
The notation is explained in \autoref{s.list.1} below.
The maximal number of conics is $432$, as in~\cite{Bauer:conics}.

\subsection{Patterns}\label{s.patterns}
As explained in \autoref{s.extensions}, a codimension~$1$ stratum in~$\famB$
is
defined by a finite index extension
%of a corank~$1$ overlattices $\lattice\(u)$ as in~\eqref{eq.Sh()},
$N\supset\lattice\(u)$, see~\eqref{eq.Sh()},
where
$u\notin\lattice$ is an extra line or
%irreducible
conic.
(We do not use the irreducibility of~$u$.)
%(Note though that the irreducibility is not assumed below.)
%(Note though that below the irreducibility is never assumed explicitly.)
If $\ls|\supp_1u|=p$ and $\ls|\supp_2u|=q$, we say that
\[*
\text{$u$ has \emph{pattern $\cn p-q$} (if $u$ is a conic) or
 \emph{$\ln p-q$} (if $u$ is a line)}.
\]

By the Hodge index theorem, $\lattice\(u)\subset\NS(X)$ must be hyperbolic.
Denoting by~$\tu$ the orthogonal projection of~$u$ to
$\lattice^\perp\otimes\Q$ and using~\eqref{eq.splitting}, we arrive at
\[*
36\tu^2=
%-(p + 2q + \Ge - 9)^2+36q+18\Ge-9,
-(p+2q+\Ge)^2+18p+72q-72,
\]
where $p$, $q$ are as above and $\Ge:=u\cdot h\in\{1,2\}$. The requirement
that $\tu^2<0$ results in \autoref{tab.Sylvester}, left (for conics) or right
(for lines). In the table,
\roster*
\item
a dot~$\cdot$ marks the pairs $(p,q)$ for which $\tu^2\ge0$,
\item
a cross~$\times$ marks the remaining
pairs ruled out by the parity condition~\eqref{eq.parity};
it is due to~\eqref{eq.parity} that we list only even values of~$p$, and
\item
a $\circ$ marks the pairs further ruled out by one of the lemmas below.
\endroster

\table
\caption{Sylvester test for conics (left) and lines (right)}\label{tab.Sylvester}
\def\bx#1{\hbox to9pt{\hss$#1$\hss}}%
\def\select#1{\if#1.\cdot\else
 \ifcase#1
  \times\or
  \bullet\or
  \circ\or
%  *\fi\fi
  \circ\fi\fi
}%
\def\r#1{\hss$\scriptstyle#1$\ \,}%
\def\q#1{\omit\bx{\scriptstyle#1}}%
\hrule height0pt
\centerline{\baselineskip10pt%
\vtop{\halign{#&&\bx{\select{#}}\cr
&\q{0}&&\q{2}&&\q{4}&&\q{6}&&\q{8}&&\q{10}&&\q{12}&&\q{14}&&\q{16}\cr
\r{0}&2&2&.&.&.&.&.&.&.&.&.&.&.&.&.&2&2\cr
\r{2}&0&.&.&.&.&.&.&.&.&.&.&.&.&.&0\cr
\r{4}&1&.&.&.&.&.&.&.&.&.&.&.&1\cr
\r{6}&1&.&.&.&.&.&.&.&.&.&1\cr
\r{8}&1&.&.&.&.&.&.&.&1\cr
\r{10}&1&2&3&.&3&2&1\cr
\r{12}&1&1&1&1&1\cr
\r{14}&0&0&0\cr
\r{16}&2\cr}}%
\quad\
\vtop{\halign{#&&\bx{\select{#}}\cr
&\q{0}&&\q{2}&&\q{4}&&\q{6}&&\q{8}&&\q{10}&&\q{12}&&\q{14}&&\q{16}\cr
&0&0&.&.&.&.&.&.&.&.&.&.&.&.&.&.&0\cr
&0&.&.&.&.&.&.&.&.&.&.&.&.&.&.\cr
&1&.&.&.&.&.&.&.&.&.&.&.&.\cr
&1&.&.&.&.&.&.&.&.&.&.\cr
&3&.&.&.&.&.&.&.&.\cr
&1&.&.&.&.&.&.\cr
&1&3&3&3&3\cr
&0&0&0\cr
&0\cr}}}
\endtable

The next two lemmas state, in particular,
that the lattice $\lattice\(u)$ depends on the
\emph{size} $\ls|\supp_2u|$ and \emph{equivalence class} $\cl{\supp_1u}$
only.
%Unlike most other statements below, Lemmas~\ref{lem.supp2}
%and~\ref{lem.class}\iref{i.equiv} hold for any polarization $h^2\in4\Z$.

\lemma\label{lem.supp2}
Let $u\notin\lattice$ be an extra conic \rom(line\rom), $\gu:=\supp_1u$,
and $\ggu:=\supp_2u$. Then, for any
set $\gv'\subset\bar\gu$, $\ls|\gv'|=\ls|\ggu|$, there is a conic
\rom(respectively, line\rom)
$v\in\lattice\(u)$  such that
$\supp_1v=\gu$ and $\supp_2v=\gv'$.
%\[*\supp_1v=\gu,\quad\supp_2v=\gv'.\]
\endlemma

\proof
The vector in question is $u+(\gs\sminus\ggu)-(\ggu\sminus\gs).$
\endproof

\lemma\label{lem.class}
Let $u\notin\lattice$ be an extra conic \rom(line\rom), $\gu:=\supp_1u$, and
$p:=\ls|\gu|$. Let, further, $\gv$ be any of the following sets\rom:
\roster
\item\label{i.equiv}
$\gv\in\cl\gu_p$, or
\item\label{i.hequiv.conic}
$\gv\in\Cl\gu_{14-p}$, if $u$ is a conic and $\supp_2u=\varnothing$, or
\item\label{i.hequiv.line}
$\gv\in\Cl\gu_{16-p}\sminus\cl\gu$, if $u$ is a line and $\supp_2u=\varnothing$.
\endroster
Then, there is a conic \rom(respectively, line\rom) $v\in\lattice\(u)$ such
that $\supp_1v=\gv$.
\endlemma

In items~\iref{i.hequiv.conic},~\iref{i.hequiv.line}, the assumption
$\supp_2u=\varnothing$ is for the sake of simplicity, as the other cases are
redundant, see Lemmas~\ref{lem.cn10-q} and~\ref{lem.ln12-q} below.

%In fact,
%the contradiction for \latin{loc\PERIOD\ cit\PERIOD} is obtained by following
%the construction in the proof of \autoref{lem.class}.

\proof[Proof of \autoref{lem.class}]
Any set~$\gv\ne\gu$ as in item~\iref{i.equiv} is of the form $\gu\sd\go$ for
some set $\go\in\CO_8$, $\ls|\gu\cap\go|=4$. Let $\go_+:=\go\cap\supp_2u$ and
pick a subset $\go_-\subset\go\cap\gu$ such that $\ls|\go_-|=\ls|\go_+|$.
Then,
$v:=u+\ds\gu\go+\go_+-\go_-$
is  as in the statement.

A set~$\gv$ as in item~\iref{i.hequiv.conic} is of the form
$\overline{\gu\sd\go}$ for some $\go\in\CK_6$ such that $\ls|\go\cap\gu|=2$,
and the conic $v:=u+\hh+\ds\gu\go$ is as required.
For item~\iref{i.hequiv.line}, consider the \BB-conic
%$\overline{\gu\sd\go}$ for some $\go\in\CK_6$ such that $\ls|\go\cap\gu|=3$.
%The \BB-conic
\[
w:=\hh-\ds\go\gu,\quad\go\in\CK_6,\ \ls|\go\cap\gu|=3;
\label{eq.BB.negative}
\]
one has $w\cdot u=-1$ and, hence, $v:=w-u$ is a line,
$\supp_1v=\overline{\gu\sd\go}$.
%Finally, any~$\gv$
%as in item~\iref{i.hequiv.line} is of the form
%$\overline{\gu\sd\go}$ for some $\go\in\CK_6$ such that $\ls|\go\cap\gu|=3$.
%The \BB-conic
%\[
%w:=\hh-\ds\go\gu
%\label{eq.BB.negative}
%\]
%has $w\cdot u=-1$, and, hence, $v:=w-u$ is a line as in the statement.
\endproof

\subsection{Extra conics}\label{s.conics.1}
Below are a few further restrictions on the supports of an
extra conic $u\notin\lattice$. By~\eqref{eq.parity},
the $1$-support $\supp_1u\in\CC_*$
is an even set.

\lemma\label{lem.char}
If $u\notin\lattice$ is an extra conic and $\gu:=\supp_1u\in\KC\cup\KK$,
then the lattice $\lattice\(u)$ has no geometric extensions.
\endlemma

\proof
Assuming the contrary, let $\ggu:=\supp_2u$ and consider the vector
\[*
\hu:=
\begin{cases}
  u-\ds\gu\varnothing+\ggu, & \mbox{if $\gu\in\KC$},\\
  \hh-u+\ds{\bar\gu}{\ggu}, & \mbox{if $\gu\in\KK$}.
\end{cases}
\]
We have $\hu\in\Sperp$ (see \autoref{s.Kummer.structure}) and, respectively,
\[*
\alignedat3
\hu^2&=\phantom{-}\tfrac12\ls|\gu|+2\ls|\ggu|-2,&
 \quad\hu\cdot h&=\phantom{0}2+\ls|\gu|+2\ls|\ggu|
 &\quad&\text{if $\gu\in\KC$},\\
\hu^2&=-\tfrac12\ls|\gu|+5,&\quad\hu\cdot h&=16-\ls|\gu|-2\ls|\ggu|
 &\quad&\text{if $\gu\in\KK$}.
\endalignedat
\]
In view of~\eqref{eq.Sperp}, the existence of $\hu\in\Sperp$ rules out
%the
patterns
$\cn p-0$, $p=0$, $6$, $8$, $16$. The few remaining cases are considered
below; we take into account~\autoref{tab.Sylvester} and use the
duality~\eqref{eq.duality} to assume that $\ls|\gu|+2\ls|\ggu|\le16$.
In all cases, $\hu^2=0$.

\emph{The pattern $\cn10-2$}:\enspace
we have $\hu\cdot h=2$, \ie, $\hu$ is $2$-isotropic,
see \autoref{def.admissible}\iref{i.isotropic}.

\emph{The patterns $\cn0-1$ and $\cn10-q$, $q=0,1$}:\enspace
by \autoref{lem.kernels}, a geometric extension of $\lattice\(u)$ must
contain  $v:=-\frac12\hu-\ds{\gs}\varnothing$ for some $\gs\in\CC_4$.
We have $\hu\cdot h=6$ or $4$; hence,
\[*
v^2=-2,\quad\text{$v\cdot h=1$ or~$2$},\quad\text{$v\cdot e=-1$ for each $e\in\gs$},
\]
resulting in a missing conic, see \autoref{def.admissible}\iref{i.missing},
or exceptional divisor $v-e$, see
\autoref{def.admissible}\iref{i.exceptional}, respectively.
\endproof

\lemma\label{lem.cn10-q}
Let $u\notin\lattice$ be an extra conic, $\gu:=\supp_1u\in\CC_{10}$, and
assume that
$\ggu:=\supp_2u\ne\varnothing$.
Then the lattice $\lattice\(u)$ has no geometric extensions.
\endlemma

\proof
Using the duality~\eqref{eq.duality}, assume that $\ls|\ggu|\le2$.
Since $\gu\notin\CK_{10}$ (see \autoref{lem.char}),
there is a set $\gs\in\CK_6$ such that
$\ls|\gu\cap\gs|=2$; by \autoref{lem.supp2}, we can also assume that
$\ggu\subset\gs$. If $\ls|\ggu|=2$, there is a $3$-element subset
$\gr\subset\gs\sminus\ggu$
such that $\ls|\gu\cap\gr|=1$;
if $\ls|\ggu|=1$, take $\gr:=\gs\sminus(\gu\cup\ggu)$.
In both cases,
we have
$u\cdot v=-1$ for the \BB-conic
\[*
v:=\hh+\ds\gs\gr\in\lattice\(u),
\]
\cf.~\eqref{eq.10-3}; hence,
$u-v$ is an exceptional divisor, see
\autoref{def.admissible}\iref{i.exceptional}.
\endproof

%\lemma\label{lem.cn12-2}
%Let $u$ be an extra conic with pattern~$\cn12-2$. Then, for any
%$\gv\in\cl{\supp_1u}_{12}$ and $\gv'\in\bar\gv$, $\ls|\gv'|=2$, there are
%at least \emph{two} conics $v_1,v_2\in\lattice\(u)$ such that
%\[*
%\supp_1v_i=\gv\quad\text{and}\quad\supp_2v_i=\gv'.
%\]
%%$\supp_1v_i=\gv$ and
%%$\supp_2v_i=\gv'$.
%\endlemma
%
%\proof
%One conic $v_1$ is given by Lemmas~\ref{lem.class}\iref{i.equiv}
%and~\ref{lem.supp2}, and the other one is
%\[*
%v_2:=h-v_1+\gv'-\bar\gv',
%\]
%\cf.~\eqref{eq.duality} and the proof of \autoref{lem.supp2}.
%\endproof

\subsection{Extra lines}\label{s.lines.1}
The remaining restrictions on
%the two supports of
an extra line $u\notin\lattice$
are given by the following lemma.
By~\eqref{eq.parity}, $\supp_1u\in\CC_*$ is an odd set.

\lemma\label{lem.ln12-q}
Let $u\notin\lattice$ be an extra line, and let $\gu:=\supp_1u$, $\ggu:=\supp_2\gu$.
If $\ls|\gu|=8$ or $\ls|\gu|=12$
and $\ggu\ne\varnothing$, then the lattice $\lattice\(u)$ is not admissible.
\endlemma

\proof
If $\ls|\gu|=8$ (and necessarily $\ggu=\varnothing$, see
\autoref{tab.Sylvester}), then $v:=2u+\ds\dm\gu$ is a $2$-isotropic vector, see
\autoref{def.admissible}\iref{i.isotropic}. Thus, assume that $\ls|\gu|=12$
and, as in the proof of \autoref{lem.class}, let $\go\in\CK_6$ be such that
$\ls|\go\cap\gu|=3$. Using \autoref{lem.supp2}, we can assume that
$\go_+:=\go\cap\supp_2u\ne\varnothing$, and pick $\go_-\subset\go\cap\gu$ so
that $\ls|\go-|=\ls|\go_+|$. It is immediate that
\[*
v:=\hh-\ds\go\gu-\go_++\go_--u
\]
is a line, but $v\cdot e=-1$ for each $e\in\go_-\ne\varnothing$, \cf.
\autoref{def.admissible}\iref{i.missing}.
\endproof

\subsection{The list of strata}\label{s.list.1}
In view of \autoref{tab.Sylvester} and Lemmas~\ref{lem.char},
\ref{lem.cn10-q}, and~\ref{lem.ln12-q}, we are left with the following eleven
patterns for an extra conic or line:
\[*
\cn4-0,\cn10-0;\cn6-0,\cn8-0;\cn12-0;\cn12-1;\cn12-2\quad\text{or}\quad
\ln4-0,\ln12-0;\ln6-0,\ln10-0.
\]
By \autoref{lem.class}\iref{i.hequiv.conic}, \iref{i.hequiv.line}, each pair of
patterns separated by a comma, \eg, $\cn4-0,\cn10-0$, results in the same
lattice $\lattice\(u)$: in the example, if $u$ has pattern~$\cn4-0$,
there exists $v\in\lattice\(u)$ with pattern $\cn10-0$ such that
$\lattice\(u)=\lattice\(v)$,
%then
%$\lattice\(u)=\lattice\(v)$ for some $v\in\lattice\(u)$ of type $\cn10-0$,
and \latin{vice versa}.

Lemmas~\ref{lem.supp2} and~\ref{lem.class}\iref{i.equiv} state
that the lattice $\lattice\(u)$ depends only on the pattern of~$u$ and
the equivalence class $\cl\gu$ of its $1$-support $\gu:=\supp_1u$.
Up to isomorphism preserving~$\dm$,
the corank~$1$ extension
$\lattice\(u)\supset\lattice$ is determined by the $\G$-orbit of~$\gu$; due to
\autoref{tab.orbits}, parity condition~\eqref{eq.parity}, and
\autoref{lem.char} eliminating the exceptional orbit $\KC\cup\KK$,
we arrive at the seven classes collected in \autoref{tab.1}.
Letting
$X\subset\famX$ be a very general member of the respective family,
listed in the table are
\roster*
\item
the notation for the codimension~$1$ stratum $\famX\subset\famB$ (for future references),
\item
the patterns of the lines and conics on $X$,
\item
the projections $\Gd_2$, $\Gd_3$ to $\discr\lattice$
(see \autoref{s.clusters} below).
\endroster
The rest of the data is common for Tables~\ref{tab.1}, \ref{tab.2},
and~\ref{tab.3}:
\roster*
\item
the numbers of lines and conics on $X$; the latter is
given by a single count if all conics are irreducible or as
$(\text{reducible})+(\text{irreducible})$ otherwise,
\item
the size of the group $G:=\Aut(\Fn X)$ of abstract graph automorphisms of the
bi-colored Fano graph $\Fn X$,
\item
the index $i_\dm:=[G:G_\dm]$ of the subgroup $G_\dm\subset G$
preserving~$\dm$ as a set,
\item
the group (reference to \autoref{tab.groups}) $G_\Go$ of symplectic
automorphisms of~$X$ and the index $[\Aut_hX:G_\Go]$ (if greater than~$1$) as
a  superscript,
\item
the determinant $\ls|\det T(X)|$ or, in \autoref{tab.3}, the lattice $T(X)$
itself,
\item
the numbers $(r,c):=\rc\famX$, see~\eqref{eq.(r,c)(Y)}, and
\smash{$(\br_i,\bc_i):=\rc\famXdm_i$}
(separately for each relative
form, \cf. \autoref{rem.forms} below),
in the form
\[
\vcenter{\halign{#,\hss\quad&if $#$\hss\cr
 both values $(r,c)\NE(\br_i,\bc_i)$&(r,c)\ne(\br_i,\bc_i),\cr
 the common value $(r,c)$&(r,c)=(\br_i,\bc_i).
\crcr}}
\label{eq.(r,c)}
\]
\endroster
A straightforward computation using \GAP~\cite{GAP4}
yields the following statement.

\proposition\label{prop.1}
For each
of the
%seven
\rom(pairs of\rom) patterns in \autoref{tab.1}, the
corank~$1$ extension
%lattice
$\lattice\(u)$ is
geometric and has no proper geometric finite index extensions. Furthermore,
the lines and conics in
the complement
$\lattice\(u)\sminus\lattice$ are precisely those
given by Lemmas~\ref{lem.supp2} and \ref{lem.class}
and duality~\eqref{eq.duality}.
%and~\ref{lem.cn12-2}.
\done
\endproposition

Concerning the last statement, in all but one cases each pair $\gu$, $\ggu$
given by the lemmas and~\eqref{eq.duality} is represented by a single
line/conic. The last case $\rstrat7$, with the pattern $\cn12-2$, is special,
as \autoref{lem.supp2} and~\eqref{eq.duality} give us two conics for each
pair.

%Note that, in
In the presence of lines
(rows~$\rstrat1$ and~$\rstrat2$ of the table),
some of the \BB-conics,
\viz. those given by~\eqref{eq.BB.negative}, become reducible.
However, the
lattice contains no new conics, reducible or not: the total conic count is
still $352$.
%; this fact is reflected in the table.

\subsection{Clusters}\label{s.clusters}
Given a geometric extension $N\supset\lattice$, we define a \emph{cluster} as
the set of all lines and conics (irreducible or reducible) in
$N\sminus\lattice$ that project to the same point in the projectivization
$\mathbb{P}\bigl((N/\lattice)\otimes\Q\bigr)$. According to \autoref{prop.1},
there are seven isomorphism classes
(henceforth referred to as the \emph{types})
of clusters listed in \autoref{tab.1} and, for each cluster~$C$,
\[*
\text{the primitive hull of $\lattice\cup C$ in~$N$ is $\lattice\(u)$ for any
$u\in C$}.
\]
%Furthermore, all
All vectors $u\in C$ share
the same (pair of) patterns and
%the same
class $\Cl{\supp_1u}$.

For an alternative description of a cluster~$C$,
observe that a vector $u\in N$ defines a
linear functional $\lattice\to\Z$, $x\mapsto x\cdot u$, and, thus, a class
in $\discr\lattice$. Since any two vectors $u,v\in C$ are related \via\
$u\pm v\in\lattice$ (\cf. Lemmas~\ref{lem.supp2} and~\ref{lem.class}), we
conclude that $C$ projects to a pair of opposite vectors
$\pm\Gd(C)\in\discr\lattice$.
%(If $C$ contains two patterns,
%we choose for $+\Gd(C)$ the image of that with shorter support.)
The components
\[*
\Gd_p:=\Gd_p(C)\in\discr_p\lattice,\quad p=2,3,
\]
of $\Gd(C)$
are characterized in \autoref{tab.1} by means of the value $\Gd_2^2\bmod2\Z$
and element $\Gd_3$ of the group
$\discr_3\lattice\cong\Z/9$ generated by $\frac19(h+\dm)$.
Note that
\[
\vcenter{\openup3pt\halign{#\hss\cr
$\Gd_2(C)$ has order~$2$ or~$4$ if $C$ consists of conics or lines,
 respectively, and\cr
$\Gd_2\ne0$ or $\hh\bmod\lattice$ (see \autoref{lem.char}).
\crcr}}
\label{eq.cluster}
\]
%\roster*
%\item
%$\Gd_2(C)$ has order~$2$ or~$4$ if $C$ consists of conics or lines,
%respectively, and
%\item
%$\Gd_2\ne0$ or $\hh\bmod\lattice$ (see \autoref{lem.char}).
%\endroster
%Observing
%%in addition
%that $\Gd_2\ne0$ or $\hh\bmod\lattice$ (see \autoref{lem.char}) and comparing
Comparing
the counts, \cf.~\eqref{eq.discr} and \autoref{tab.orbits}, we arrive at the
following statement.

\proposition\label{prop.cluster}
The common support $\Cl{\supp_1u}$, $u\in C$, of any cluster $C\subset N$ and, hence,
the sublattice $\lattice\(u)\subset N$ is determined
by the image $\Gd_2(C)$ and, in the case of types
$\rstrat3$, $\rstrat4$, $\rstrat5$, $\rstrat7$, by any one of the
following\rom:
\roster*
\item
the type specification
\rom($\rstrat3$, $\rstrat4$, $\rstrat5$, or $\rstrat7$\rom),
or, equivalently,
\item
the image $\Gd_3(C)\in\discr_3\lattice$.
\done
\endroster
\endproposition

In the next section (see \autoref{cor.clusters}) we observe that
distinct clusters $C_1$, $C_2$
in any \emph{geometric} extension $N\supset\lattice$ have
distinct images $\pm\Gd_2(C_1)\ne\pm\Gd_2(C_2)$.

\section{Strata of codimension $2$}\label{S.codim2}

\table
\caption{Strata of codimension~$2$ \rm(see~\autoref{s.tab2} and \autoref{s.list.1})}\label{tab.2}
\centerline{\vbox{\small\offinterlineskip
\tabdefs
\def\SAME{2}%
\tabskip0pt plus 1fil
\halign{\strut\quad$#$\hss\sep\tabskip0pt&
 \strut\sep\sep\hss$#$\hss\sep&
 \strut\sep\hss$#$\hss\sep&
 \strut\sep\hss$#$\hss\sep&
 \strut\sep\hss$#$\hss\sep&
 \strut\sep\hss$#$\hss\sep&
 \strut\sep\hss$#$\hss\sep&
 \strut\sep\hss$#$\hss\sep&
 \strut\sep$#$\hss\sep\tabskip0pt plus 1fil\cr
\noalign{\hrule\kern3pt}
\text{Clusters}&\circ&\text{Lines}&\text{Conics}&\ls|G|&i_\dm&
 G_\Go&\ls|{\det}|&\omit\hss$(r,c)$\hss\cr
\noalign{\kern1pt\hrule\kern3pt}
\s1,\s1,\s4&1,2&16&48\+384&512&4&\(16,14)\^2&224&(1,0)\EQ(1,0)\cr
\s1,\s1,\s7&3,2&16&48\+352&256&2&\(16,14)\^2&256&(1,0)\EQ(1,0)\cr
\s1,\s2,\s3,\s6&3,0,2&40&240\+184&192&2&\(16,14)\^1&144&(1,0)\EQ(1,0)\cr
\noalign{\kern1pt\hrule\kern3pt}
\s1,\s3&0&8&16\+376&384&4&\(16,14)\^1&336&(1,0)\EQ(1,0)\cr
\s1,\s3&2&8&16\+376&256&4&\(16,14)\^1&384&(1,0)\EQ(1,0)\cr
\s1,\s4&0&8&16\+416&384&4&\(16,14)\^2&304&(1,0)\EQ(1,0)\cr
\s1,\s5&0&8&16\+344&192&2&\(16,14)\^1&464&(1,0)\EQ(1,0)\cr
\s1,\s5&0&8&16\+344&192&2&\(16,14)\^1&400&(1,0)\EQ(1,0)\cr
\s1,\s5&2&8&16\+344&128&2&\(16,14)\^1&576&(1,0)\EQ(1,0)\cr
\s1,\s6&2&8&16\+368&384&4&\(16,14)\^1&400&(1,0)\NE(2,0)\cr
\s1,\s6&0&8&16\+368&256&4&\(16,14)\^1&416&(1,0)\EQ(1,0)\cr
\s1,\s7&0&8&16\+384&384&2&\(16,14)\^2&400&(1,0)\EQ(1,0)\cr
\s2,\s3&2&32&160\+232&768&4&\(16,14)\^1&192&(1,0)\EQ(1,0)\cr
\s2,\s4&2&32&160\+272&768&4&\(16,14)\^2&160&(1,0)\EQ(1,0)\cr
\s2,\s5&0&32&160\+200&192&2&\(16,14)\^1&272&(1,0)\EQ(1,0)\cr
\s2,\s5&2&32&160\+200&384&2&\(16,14)\^1&256&(1,0)\EQ(1,0)\cr
\s2,\s6&0&32&160\+224&768&4&\(16,14)\^1&224&(1,0)\EQ(1,0)\cr
\s2,\s7&0&32&160\+240&384&2&\(16,14)\^2&208&(1,0)\EQ(1,0)\cr
\noalign{\kern1pt\hrule\kern3pt}
\same{3} \s3,\s3&0&&432&512&2&\(16,14)\^1&576&(1,0)\EQ(1,0)\cr
\s3,\s3,\s3&2&&472&4608&8&\(16,14)\^1&432&(1,0)\NE(2,0)\cr
\same{1} \s3,\s3,\s7&2&&480&384&2&\(16,14)\^2&432&(1,0)\EQ(1,0)\cr
\same{1} \s3,\s4,\s5&2&&480&384&4&\(16,14)\^2&432&(1,0)\EQ(1,0)\cr
\s3,\s4,\s6&0&&504&1024&8&\(16,14)\^2&384&(1,0)\EQ(1,0)\cr
\s3,\s5&0&&400&512&4&\(16,14)\^1&768&(1,0)\EQ(1,0)\cr
\s3,\s5&2&&400&384&4&\(16,14)\^1&816&(1,0)\EQ(1,0)\cr
\same{3} \s3,\s5,\s6&0&&432&512&4&\(16,14)\^1&576&(1,0)\EQ(1,0)\cr
\same{4} \s3,\s6&2&&424&384&4&\(16,14)\^1&624&(1,0)\EQ(1,0)\cr
\s3,\s7&0&&440&1024&4&\(16,14)\^2&576&(1,0)\EQ(1,0)\cr
\s4,\s4,\s6&0&&544&2048&8&\(32,27)\^2&320&(1,0)\EQ(1,0)\cr
\maxr1 \s4,\s4,\s7&2&&560&1152&6&\(48,50)\^2&304&(1,0)\EQ(1,0)\cr
\s4,\s5&0&&440&512&4&\(16,14)\^2&640&(1,0)\EQ(1,0)\cr
\s4,\s5,\s5&2&&448&768&4&\(16,14)\^2&560&(1,0)\EQ(1,0)\cr
\same{2} \s4,\s6&2&&464&384&4&\(16,14)\^2&496&(1,0)\EQ(1,0)\cr
\s4,\s7&0&&480&1024&4&\(32,27)\^2&448&(1,0)\EQ(1,0)\cr
\s5,\s5&0&&368&512&2&\(16,14)\^1&1088&(1,0)\EQ(1,0)\cr
\s5,\s5&0&&368&512&2&\(16,14)\^1&1024&(1,0)\EQ(1,0)\cr
\s5,\s5&2&&368&384&2&\(16,14)\^1&1200&(1,0)\EQ(1,0)\cr
\s5,\s5,\s7&2&&416&384&2&\(16,14)\^2&688&(1,0)\EQ(1,0)\cr
\s5,\s6&0&&392&512&4&\(16,14)\^1&896&(1,0)\EQ(1,0)\cr
\s5,\s6&2&&392&384&4&\(16,14)\^1&880&(1,0)\EQ(1,0)\cr
\same{4} \s5,\s6,\s6&2&&424&384&2&\(16,14)\^1&624&(1,0)\EQ(1,0)\cr
\s5,\s7&0&&408&512&2&\(16,14)\^2&832&(1,0)\EQ(1,0)\cr
\s6,\s6&0&&416&1536&6&\(16,14)\^1&704&(1,0)\EQ(1,0)\cr
\same{2} \s6,\s6,\s7&2&&464&384&2&\(16,14)\^2&496&(1,0)\EQ(1,0)\cr
\s6,\s7&0&&432&1024&4&\(16,14)\^2&640&(1,0)\EQ(1,0)\cr
\s7,\s7&0&&448&2048&2&\(32,27)\^2&576&(1,0)\EQ(1,0)\cr
\s7,\s7,\s7&2&&496&2304&2&\(48,50)\^2&432&(1,0)\EQ(1,0)\cr
\noalign{\kern1pt\hrule}
\crcr}}}
\endtable

Next step is the analysis of the double (self-)intersections of the seven
strata of codimension~$1$ described in \autoref{tab.1}.
%This is a
%straightforward computation (\cf. \autoref{s.proof2} below), and the
%resulting $47$ relative forms of the $43$
%$1$-parameter (up to projective transformation) equiconical families of
%Barth--Bauer quartics are listed in \autoref{tab.2}.
The resulting $47$ relative forms of the $43$ equiconical families are shown
in \autoref{tab.2}. The maximal number of conics is $560$; it is in the row
marked with a \maxreal1 (see also \autoref{th.real}).

To compile the table, we analyze all corank~$2$ extensions
\[*
\lattice\(u_1,u_2),\quad u_i\in C_i,\ i=1,2,
\]
where $C_1\ne C_2$ is a pair of distinct clusters (see \autoref{s.clusters}),
which are represented by their classes $\Ga_i:=\Gd_2(C_i)\in\discr_2\lattice$, $i=1,2$.
The candidates are found as follows:
\roster*
\item
pick one of the $27$ $\G$-orbits of pairs $(\Ga_1,\Ga_2)$
satisfying~\eqref{eq.cluster},
\item
for each $\Ga_i$ with $\Ga_i^2=1\bmod2$, pick a type, see
\autoref{prop.cluster}, and
\item
specify the product $u_1\cdot u_2$, see item~\iref{ex.Gram} in
\autoref{s.extensions}.
\endroster
For the latter, and for the construction of~$\lattice\(u_1,u_2)$, we need to
select a pair of representatives $u_i\in C_i$. The values for $u_1\cdot u_2$
are limited by \autoref{lem.convention} and the Hodge index theorem: the
orthogonal
projection of
%$(u_1,u_2)$
$\Z u_1+\Z u_2$ to $\lattice^\perp\otimes\Q$
%is
must be negative definite.

The resulting list is analyzed as explained in \autoref{s.remarks};
%All forbidden pairs are ruled out by simple arguments similar to
%\autoref{S.codim1}, but
we omit the details.
%but, in view of the large number of cases to be
%considered, we omit the details.

\subsection{Notation in \autoref{tab.2}}\label{s.tab2}
We use the same notation as in \autoref{tab.1} (see \autoref{s.list.1}),
%except that
but instead of the
first four columns we list the types (as references to \autoref{tab.1})
of the clusters~$C_i$ found in~$N$ and the products~$\circ$ of the
images $\Ga_i:=\Gd_2(C_i)\in\discr_2\lattice$.
%of the generating clusters.
Each product is $n/4\bmod\Z$, where $n$ is the number given in the table.
%Besides, the following convention is used:
We show:
\roster*
%\item
%for two elements $\Ga_1,\Ga_2$, the product $\Ga_1\cdot\Ga_2$;
%%if there are two elements $\Gd'$, $\Gd''$, we indicate
%%the product $\Gd'\cdot\Gd''$;
\item
for three \emph{order~$2$} elements $\Ga_1,\Ga_2,\Ga_3$,
the common product $\Ga_i\cdot\Ga_j$, $i\ne j$;
%if there are three elements $\Gd'$, $\Gd''$, $\Gd'''$ of order~$2$
%(conics), the
%number shown in the table refers to the common product
%$\Gd'\cdot\Gd''=\Gd'\cdot\Gd'''=\Gd''\cdot\Gd'''$;
\item
in all other cases, the products $\Ga_1\cdot\Ga_i$, $i\ge2$.
%in the presence of two order~$4$ elements $\Ga_1,\Ga_2$, just the product
%$\Ga_1\cdot\Ga_2$.
%if there are two elements $\Gd'$, $\Gd''$ of order~$4$ (lines) and one or
%several elements $\Ga,\ldots$ of order~$2$ (conics), the numbers
%$m;n,\ldots$ given in the table refer to
%$\Gd'\cdot\Gd''$ and $\Gd'\cdot\Ga=\Gd''\cdot\Ga,\ldots$, respectively.
\endroster
%With the extra observation that
%the images $\Gd_2(C)$ of all clusters $C\subset N$
%generate a subgroup $\Z/2\times\Z/2\not\ni\hh\bmod\lattice$
%or $\Z/4\times\Z/2$, see
In view of \autoref{cor.2-subgroup} below,
these data determine the collection of elements up to the
induced action of~$\G$.

For strata of codimension~$2$,
an abstract Fano
graph~$\graph$ may admit several relative forms.
They are listed in separate rows and prefixed with
equal superscripts: \eg, the two rows prefixed with~\rsame21
represent the same abstract Fano graph.

The row prefixed with a \maxreal1 contains the family maximizing the number of
real conics in a real Barth--Bauer quartic
(see \autoref{th.real}).

\subsection{A few consequences of \autoref{tab.2}}\label{s.consequences.2}
In spite of the large number of classes, one can still observe a few
common properties.
We omit their proofs, as they mainly consist of ruling out a
large number of simple
cases similar to \autoref{S.codim1}.

\corollary\label{cor.clusters}
In any \emph{geometric} extension $N\supset\lattice$,
distinct clusters $C_1,C_2\subset N$
project to distinct pairs
$\{\pm\Gd_2(C_1)\}\ne\{\pm\Gd_2(C_2)\}\subset\discr_2\lattice$.
\done
\endcorollary

\corollary\label{cor.2-subgroup}
Given a geometric corank~$2$ extension $N\supset\lattice$
generated over~$\Q$ by lines and conics,
the images
$\Gd_2(C)$ of all clusters $C\subset N$ generate a subgroup
\[*
\Z/2\times\Z/2\not\ni\hh\bmod\lattice\quad\text{or}\quad\Z/4\times\Z/2
\]
in $\discr_2\lattice$.
\done
\endcorollary

\corollary\label{cor.2-generated}
A geometric corank~$2$ extension $N\supset\lattice$ that is
generated over~$\Q$ by lines and conics has two to four clusters.
Given
any two clusters $C_1\ne C_2\subset N$
\rom(except $\rstrat3,\rstrat6$ in the third row of \autoref{tab.2}\rom)
and any pair of representatives $u_i\in C_i$,
$i=1,2$, one has $N=\lattice\(u_1,u_2)$.
%Hence, $N$ is generated by lines and conics over~$\Z$.
\done
\endcorollary

Immediately by the definition of clusters, the excessive (\ie, beyond the
common $352$ conics) line/conic counts of~$N$ are the sums of those over all
clusters $C\subset N$. As in the case of codimension~$1$, in the presence of
lines (clusters~$\rstrat1$ and~$\rstrat2$), some of the \BB-conics become
reducible.
%, and the number of those is also additive over all clusters (\cf.
%\autoref{tab.1}).
However, starting from codimension~$2$, some of the extra
conics may also be reducible; for example, in the third row of
\autoref{tab.2}, \emph{each} of the $32$ conics in the last cluster~$\rstrat6$ is
reducible.
Since each reducible conic lies in a corank~$2$ extension, we can state the
following general count:
\[
\aligned
\strut\#\{\text{reducible conics}\}
&=\text{$16$ \BB-conics for each cluster~$\rstrat1$}\\
&\mathrel+\text{$160$ \BB-conics for each cluster~$\rstrat2$}\\
&\mathrel+\text{$16$ type~$\rstrat4$ or~$\rstrat7$ conics for each pair~$\rstrat1,\rstrat1$}\\
&\mathrel+\text{$32$ type~$\rstrat3$}
 +\text{$32$ type~$\rstrat6$ conics for each pair~$\rstrat1,\rstrat2$}.
\endaligned
\]

\section{Rigid quartics (codimension $3$)}\label{S.codim3}

\table
\caption{Rigid quartic with many conics \rm(see~\autoref{s.list.3} and \autoref{s.list.1})}\label{tab.3}
\centerline{\vbox{\small\offinterlineskip
\tabdefs
\def\SAME{3}%
\tabskip0pt plus 1fil
\halign{\strut\quad$#$\hss\sep\tabskip0pt&
 \strut\sep\hss$#$\hss\sep&
 \strut\sep\hss$#$\hss\sep&
 \strut\sep\hss$#$\hss\sep&
 \strut\sep\hss$#$\hss\sep&
 \strut\sep\hss$#$\hss\sep&
 \strut\sep\hss$#$\hss\sep&
 \strut\sep$#$\hss\sep\tabskip0pt plus 1fil\cr
\noalign{\hrule\kern3pt}
\text{Clusters}&\text{Lines}&\text{Conics}&\ls|G|&i_\dm&
 G_\Go&T&\omit\hss$(c,r)$\hss\cr
\noalign{\kern1pt\hrule\kern3pt}
\s4,\s4,\s4,\s4,\s6,\s7,\s7&&800&15360&60&\(960,11357)\^2&[4,0,40]&(1,0)\EQ(1,0)\cr
\s4,\s4,\s4,\s7,\s7,\s7&&736&3072&16&\(192,1023)\^2&[12,0,16]&(1,0)\NE(0,1)\cr
\s4,\s4,\s4,\s6,\s6,\s7&&704&768&12&\(96,227)\^2&[8,4,28]&(0,1)\EQ(0,1)\cr
\s3,\s3,\s4,\s4,\s5,\s6,\s7&&680&768&12&\(48,50)\^2&[12,4,20]&(0,1)\NE(0,2)\cr
\s3,\s4,\s4,\s6,\s6,\s7&&664&768&12&\(48,50)\^2&[4,0,60]&(1,0)\NE(2,0)\cr
\s7,\s7,\s7,\s7,\s7,\s7&&640&3072&2&\(192,1023)\^2&[4,0,72]&(1,0)\EQ(1,0)\cr
\same{1} \s3,\s3,\s4,\s6,\s6,\s7&&624&256&4&\(32,27)\^2&[12,4,24]&(0,2)\EQ(0,2)\cr
\same{1} \s3,\s4,\s4,\s5,\s6,\s6&&624&256&8&\(32,27)\^2&[12,4,24]&(0,2)\NE(0,4)\cr
\same{2} \s3,\s3,\s3,\s7,\s7,\s7&&616&384&2&\(48,50)\^2&[4,0,72]&(2,0)\EQ(2,0)\cr
\same{2} \s3,\s4,\s4,\s5,\s5,\s7&&616&384&6&\(48,50)\^2&[4,0,72]&(2,0)\EQ(2,0)\cr
\maxr2 \s4,\s4,\s7,\s7&&608&768&6&\(96,227)\^2&[4,0,76]&(1,0)\EQ(1,0)\cr
&&&&&&[16,4,20]&(0,1)\EQ(0,1)\cr
\same{3} \s3,\s4,\s4,\s5,\s5,\s5,\s6&&608&512&8&\(32,27)\^2&[8,0,36]&(1,0)\NE(2,0)\cr
\same{3} \s3,\s3,\s3,\s3,\s7,\s7&&608&512&2&\(32,27)\^2&[8,0,36]&(1,0)\EQ(1,0)\cr
\same{4} \s3,\s3,\s3,\s4,\s5,\s5,\s6&&600&512&8&\(16,14)\^2&[12,0,24]&(1,0)\NE(4,0)\cr
\same{4} \s3,\s3,\s3,\s3,\s3,\s7&&600&512&4&\(16,14)\^2&[12,0,24]&(1,0)\NE(0,1)\cr
\same{5} \s4,\s4,\s6,\s7&&592&384&6&\(48,50)\^2&[12,0,28]&(0,2)\EQ(0,2)\cr
\same{5} \s6,\s6,\s6,\s7,\s7,\s7&&592&384&2&\(48,50)\^2&[12,0,28]&(0,2)\EQ(0,2)\cr
\s4,\s4,\s5,\s5,\s5,\s7&&584&384&6&\(48,50)\^2&[8,0,44]&(0,1)\EQ(0,1)\cr
\s4,\s7,\s7,\s7&&576&1536&4&\(96,227)\^2&[12,0,28]&(2,0)\EQ(2,0)\cr
\s3,\s3,\s4,\s5,\s5,\s7&&576&512&8&\(32,27)\^2&[12,4,28]&(0,1)\NE(0,2)\cr
\same{6} \s4,\s4,\s6,\s6&&576&512&8&\(32,27)\^2&[4,0,88]&(1,0)\NE(0,1)\cr
\same{6} \s6,\s6,\s6,\s6,\s7,\s7&&576&512&2&\(32,27)\^2&[4,0,88]&(1,0)\EQ(1,0)\cr
\s1,\s1,\s2,\s3,\s3,\s4,\s6,\s6&48&336\+240&256&4&\(16,14)\^2&[4,0,20]&(1,0)\NE(0,1)\cr
\s4,\s4,\s5,\s7&&568&384&6&\(48,50)\^2&[4,0,108]&(1,0)\EQ(1,0)\cr
&&&&&&[16,4,28]&(0,1)\EQ(0,1)\cr
\s1,\s4,\s4,\s7&8&16\+544&576&6&\(48,50)\^2&[4,2,52]&(2,0)\EQ(2,0)\cr
\s1,\s1,\s1,\s4,\s4,\s7&24&96\+464&192&6&\(48,50)\^2&[4,2,32]&(1,0)\EQ(1,0)\cr
&&&&&&[8,2,16]&(0,1)\EQ(0,1)\cr
\s2,\s4,\s4,\s7&32&160\+400&576&6&\(48,50)\^2&[4,2,28]&(1,0)\EQ(1,0)\cr
\same{7} \s3,\s3,\s5,\s6,\s6,\s7&&552&128&2&\(16,14)\^2&[4,0,92]&(1,0)\EQ(1,0)\cr
&&&&&&[12,4,32]&(0,1)\EQ(0,1)\cr
\same{7} \s3,\s4,\s5,\s5,\s6,\s6&&552&128&4&\(16,14)\^2&[4,0,92]&(1,0)\NE(0,1)\cr
&&&&&&[12,4,32]&(0,1)\NE(0,2)\cr
\same{7} \s3,\s3,\s4,\s5,\s6&&552&128&4&\(16,14)\^2&[4,0,92]&(1,0)\NE(0,1)\cr
&&&&&&[12,4,32]&(0,1)\NE(0,2)\cr
\s4,\s6,\s6,\s7&&544&512&8&\(32,27)\^2&[16,8,28]&(0,1)\NE(0,2)\cr
\s1,\s4,\s4,\s6&8&16\+528&512&8&\(32,27)\^2&[8,4,28]&(1,0)\NE(2,0)\cr
\s1,\s1,\s4,\s4,\s6&16&48\+496&512&8&\(32,27)\^2&[4,0,40]&(1,0)\NE(0,1)\cr
\s1,\s1,\s2,\s3,\s3,\s6,\s6,\s7&48&336\+208&64&2&\(16,14)\^2&[4,2,24]&(1,0)\EQ(1,0)\cr
&&&&&&[8,2,12]&(0,1)\EQ(0,1)\cr
\noalign{\kern1pt\hrule}
\crcr}}}
\endtable

There remains to analyze the ``triple intersections'' of the strata, \ie,
projectively rigid (\latin{aka} singular) quartics. The computation
%outlined in \autoref{s.proof3} below
results in $285$ relative forms of $211$ abstract
isomorphism classes of Fano graphs, most represented by
several projective equivalence classes: up to projective equivalence, there
are
\roster*
\item
$208$ real and $189$ pairs of complex conjugate quartics~$X$, and
\item
$231$ real and $682$ pairs of complex conjugate pairs $(X,\dm)$.
\endroster
The converse also holds: a singular $K3$-surface may admit several
embeddings to~$\Cp3$ as a Barth--Bauer quartic, sometimes with non-isomorphic
Fano graphs (see, \eg, the two occurrences of the
transcendental lattice $T=[4,0,40]$
in \autoref{tab.3}).

The relative forms with more than $536$ conics are shown in \autoref{tab.3},
sorted by the decreasing of the total conic count;
a complete list is found
in~\cite{degt:4Kummer.tables}.
%%\[*
%%\text{\href{http://www.fen.bilkent.edu.tr/~degt/papers/4Kummer_tab.pdf}
%% {\tt http://www.fen.bilkent.edu.tr/\char`\~degt/papers/4Kummer\_tab.pdf}}.
%%\]
%The maximal number of lines in a Barth--Bauer quartic is~$48$; both line
%maximizing Fano graphs are presented in \autoref{tab.3}.

The computation is similar to \autoref{S.codim2}:
\roster*
\item
pick one of the $134$ $\G$-orbits of triples $(\Ga_1,\Ga_2,\Ga_3)$ such that
each pair is as in \autoref{cor.2-subgroup}; in particular, all $\Ga_i$ are
distinct and satisfy~\eqref{eq.cluster},
\item
for each $\Ga_i$ with $\Ga_i^2=1\bmod2$, pick a type, see
\autoref{prop.cluster}, and
\item
specify the product $u_i\cdot u_j$, $i\ne j$, see item~\iref{ex.Gram} in
\autoref{s.extensions}.
\endroster
In the last two steps, each decorated pair $(\Ga_i,\Ga_j)$ must be one of those
in \autoref{tab.2}.

Then, for each decorated triple $(\Ga_1,\Ga_2,\Ga_3)$, we check if the
corresponding lattice $\lattice\(u_1,u_2,u_3)$ is hyperbolic, \cf.
\autoref{S.codim2}, and, if it is, analyze it as explained in
\autoref{s.remarks}.

\subsection{Notation in \autoref{tab.3}}\label{s.list.3}
The notation is similar to that of \autoref{tab.2}, see
\autoref{s.tab2} and~\autoref{s.list.1}, but we omit the
products~$\circ$ in $\discr_2\lattice$
%(which are no longer constant)
and show the transcendental lattices $T(X)$,
see~\eqref{eq.T}, instead of $\ls|\det T|$.
Each lattice,
together with its respective counts $(r,c)$ (itemized
%further
according to
the type of~$T$),
occupies a separate row.

Marked with a \maxreal2 is the quartic with the maximal number $560$ of real
conics.

\remark\label{rem.forms}
Recall that the numbers $\rc\famX$ and \smash{$\rc\famXdm_i$},
see~\eqref{eq.(r,c)(Y)}, may differ;
they are listed
in the last column
as explained in~\eqref{eq.(r,c)}, and in \autoref{tab.3} they
are further itemized according to the isomorphism type of the transcendental
lattice.

It should be understood that
the absolute numbers $\rc\famX$
are the (first pair of) values $(r,c)$
\emph{common to all rows} related to~$\graph$, whereas
%the relatives ones
\smash{$\rc\famXdm_i$}
are the (second pair of) values $(\br_i,\bc_i)$ listed
separately in the row corresponding to the
%particular
relative form
$(\graph,\dm_i)$. These latter pairs are to be summed up to obtain the
total numbers \smash{$\rc\famXdm(\graph)$}. This convention
is illustrated in the next example (see also \autoref{cor.real} below).
\endremark

\example\label{ex.forms}
We illustrate \autoref{rem.forms} by a few examples, referring to
\autoref{tab.3}. Since the quartics
%in \autoref{tab.3}
are rigid, we
%actually
speak about the
%number of
(projective equivalence classes of)
quartics rather than components.
Denote by $\BS(X)$ the set of subgraphs $\dm\subset\Fn X$, so that the
relative
forms are the $\Aut(\Fn X)$-orbits on $\BS(X)$.

In the first row, the two pairs are equal:
\smash{$\rc\famX=\rc\famXdm=(1,0)$}; thus, there is a single quartic~$X$
with this Fano graph~$\graph$, the group $\Aut_hX$ acts transitively on
%the subgraphs $\dm\subset\graph$,
$\BS(X)$,
and $X$ admits a real structure
preserving a point $\dm\in\BS(X)$.

In the second row, $\rc\famX=(1,0)$, whereas \smash{$\rc\famXdm=(0,1)$}.
Hence, there is a single quartic~$X$ with transitive $\Aut_hX$-action on
$\BS(X)$, this quartic is real, but there is no real
structure preserving a point $\dm\in\BS(X)$.

The two rows prefixed with \rsame31 represent two real forms of the same
graph:
\[*
\rc\famX=(0,2);\quad
\smash{\rc\famXdm_1=(0,2)},\
\smash{\rc\famXdm_2=(0,4)}\ \Longrightarrow\
\smash{\rc\famXdm=(0,6)}.
\]
(From this point on, the verbal description is left to the reader.)

For the two rows prefixed with \rsame32,
$\rc\famX=\smash{\rc\famXdm_1}=\smash{\rc\famXdm_2}=(2,0)$,
\smash{$\rc\famXdm=(4,0)$}.

Finally, for the three rows prefixed with~\rsame37, we have
$\rc\famX=(1,1)$: there is one real quartic, with $T\cong[4,0,92]$, and a
pair of conjugate ones, with $T\cong[12,4,32]$.
(Most likely, the quartics are Galois conjugate over an algebraic number
field, but I would refrain from an affirmative statement.)
The other counts
\[*
\smash{\rc\famXdm_1=(1,1)},\
\smash{\rc\famXdm_2=\rc\famXdm_3=(0,3)}\ \Longrightarrow\
\smash{\rc\famXdm=(1,7)}
\]
can also be itemized further according to the two types of~$T(X)$.
\endexample

\subsection{Consequences of \autoref{tab.3}\pdfstr{}{
 \rm(see also~\cite{degt:4Kummer.tables})}}\label{s.consequences.3}
Unfortunately, there seems to be no simple way to interpret \autoref{tab.3} in
terms of \autoref{tab.2}.
For example, some triples $C_1$, $C_2$, $C_3$ of
clusters independent over~$\lattice$ may project to dependent elements
$\Gd_2(C_i)\in\discr_2\lattice$,
generating a subgroup of length~$2$. For this reason, for a
typical
corank~$3$ extension $N\supset\lattice$,
%it is \emph{not} true
one cannot assert
that any decorated
(by the types of the clusters) length~$2$
subgroup in the image of~$\Gd_2$ is among the groups listed in \autoref{tab.2}.
Nor can one assert that $N=\lattice\(u_1,u_2,u_3)$ for representatives
$u_i\in C_i$ of
\emph{any} triple of independent clusters (\cf. \autoref{cor.2-generated}):
some of these lattices do have proper geometric finite index
extensions.
%generates the
%extension over~$\Z$.

Still, one can state certain limited analogues of
Corollaries~\ref{cor.2-subgroup} and~\ref{cor.2-generated}.

\corollary\label{cor.3-subgroup}
Given a geometric corank~$3$ extension $N\supset\lattice$
generated over~$\Q$ by lines and conics,
the images
$\Gd_2(C)$ of all clusters $C\subset N$ generate a subgroup
\[*
\Z/2\times\Z/2\times\Z/2\not\ni\hh\bmod\lattice
\quad\text{or}\quad\Z/4\times\Z/2\times\Z/2
\]
in $\discr_2\lattice$.
\done
\endcorollary

\corollary\label{cor.3-generated}
If a geometric corank~$3$ extension $N\supset\lattice$ is
generated by lines and conics over~$\Q$, then so it is over~$\Z$.
\done
\endcorollary

%\latin{A priori},
By definition,
the fact that a connected component
$\famR\subset\famX$ or \smash{$\famXdm$}
is real means that
%a representative
any $X\in\famR$ is
projectively equivalent to~$\bar X$,
\ie, $X$ has an anti-holomorphic automorphism, but the latter does not need
to be involutive.
%; in other words, there is an
%anti-holomorphic automorphism $\rs\:X\to X$, but this $\rs$ does not need to
%be involutive.
%, but this automorphism does not need to be involutive.
(If $\codim\famX\le2$, even less can be
stated: $X$ and $\bar X$ are related by an equiconical deformation;
\cf. \cite{degt:geography,Cisem:real} for examples of
real
strata without real
points in other similar problems.)
%.)
Nevertheless, a case-by-case analysis shows that real components of
codimension~$3$ do consist of real quartics~$X$ or, respectively, pairs
$(X,\dm)$.

\corollary\label{cor.real}
Each real component $\famR\subset\famX$ \rom(resp.\
\smash{$\famR\subset\famXdm$}\rom)
of codimension~$3$ in $\famB$
consists of real quartics $X$
\rom(resp.\ pairs $(X,\dm)$\rom).
\done
\endcorollary

\section{Missing proofs}\label{S.proofs}

In this concluding section, we fill in a few missing details to complete the
proofs of the principal results stated in \autoref{s.results}.

\subsection{Proof of \autoref{th.main}}\label{proof.main}
The only statement that needs proof is the fact that the unique Barth--Bauer
quartic~$X$ with $800$ conics (see the first row in \autoref{tab.3}) is given by
Mukai's polynomial as in the theorem.
(This was observed by X.~Roulleau, private communication.)
From the computation, we know that
$X$ admits a faithful symplectic action of Mukai's group $M_{20}$.

Let $\famM$ be the $1$-parameter family of (generically non-algebraic)
$K3$-surfaces with a faithful symplectic action of $G:=M_{20}$.
All actions have isomorphic covariant lattices $(\Fix G_*)^\perp$ (see, \eg,
\cite{Hashimoto}; here, $G_*$ is the induced action on $H_2$), and this
common lattice~$S$ is
found as
the orthogonal complement
$h^\perp\subset\NS(X)$.
We have
\[*
\discr_2S
 \cong\bmatrix1&\frac12\\\frac12&1\endbmatrix\oplus\bigl[\tfrac18\bigr],
\qquad
\discr_5S\cong\bigl[\tfrac25\bigr],
\]
so that $\NS(X)$ is the index~$2$ extension of $S\oplus\Z h$ \via\
$\tilde\Ga_1+\hh$
for a certain vector
$\tilde\Ga_1\in S\dual$ such that $\tilde\Ga_1^2=1\bmod2\Z$.

\lemma\label{lem.S.M20}
The natural homomorphism
$\dd_S\:\!\OG(S)\onto\Aut(\discr S)$ is surjective.
\endlemma

\proof
In the notation above, the group $\Aut(\Fn X)$ restricts to~$S$ and its image
under~$\dd_S$
is an index~$6$ maximal subgroup of $\Aut(\discr S)$ which, expectedly, fixes
the order~$2$ element
\[*
\Ga_1:=(\tilde\Ga_1\bmod S)\in\discr_2S.
\]
Therefore, it suffices to show that the full group $\OG(S)$ acts transitively
on the set of all six elements $\Ga_1,\ldots,\Ga_6\in\discr_2S$ of square
$1\bmod2\Z$.
To this end, pick $\tilde\Ga_i\in S\dual$ such that
$\Ga_i=\tilde\Ga_i\bmod S$ and consider the index~$2$ extension
$N_i\supset S\oplus\Z h$, $h^2=4$, generated by $\tilde\Ga_i+\hh$.
A computation shows that each $N_i\ni h$ is geometric and its Fano graph
$\Fn(N_i,h)$ consists of $800$ conics and contains a copy of~$\dm$.
Thus, due to the uniqueness given by \autoref{tab.3}, there is a lattice
isomorphism
\[*
N_i\overset\cong\longto N_1=\NS(X)
\]
preserving~$h$, \ie, an element of $\OG(S)$ taking $\Ga_i$ to $\Ga_1$.
\endproof

Given a primitive isometry $\iota\:S\into\L$, the genus of $S^\perp$
(\cf. \autoref{prop.Nikulin})
consists
of a
single isomorphism class:
\[*
S^\perp\cong
\bmatrix4&0&2\\
0&4&2\\
2&2&12\endbmatrix.
\]
Hence, by \autoref{prop.Nikulin} and \autoref{lem.S.M20}, the isometry is
unique up to $\OG(S)$-equivalence (\cf. Hashimoto~\cite{Hashimoto}).
The next corollary is proved in Bonnaf\'{e}--Sarti~\cite{Bonnafe.Sarti}
under somewhat stronger hypotheses,
assuming a faithful action of the full group $\Aut_hX_4$
(which is an index~$2$ extension of $M_{20}$, see \autoref{tab.3}).

\corollary
%[\cf. Bonnaf\'{e}--Sarti~\cite{Bonnafe.Sarti}, where stronger hypotheses are used]
\label{cor.X4}
Up to projective transformation, there is a unique quartic surface
$X_4\subset\Cp3$ admitting a faithful symplectic action of the Mukai
group~$M_{20}$.
\endcorollary

\proof
It suffices to observe that $\OGplus(S^\perp)$ acts transitively on
the set of
the four
square~$4$ vectors in $S^\perp$ and use an obvious analogue of
\autoref{prop.absolute}.
\endproof

It follows that any quartic in~$\Cp3$
admitting a faithful symplectic action of
$M_{20}$ is projectively equivalent to the example found in
Mukai~\cite{Mukai}. This observation completes the proof of
\autoref{th.main}.
\qed

\medskip
For future references, we state an analogue of \autoref{cor.X4} for
octics.
For the proof, we merely observe that the group
$\OGplus(S^\perp)$ also acts transitively on
the set of
the four
square~$8$ vectors in $S^\perp$.

\corollary\label{cor.X8}
Up to projective transformation, there is a unique
octic $K3$-surface
$X_8\subset\Cp5$ admitting a faithful symplectic action of the Mukai
group~$M_{20}$.
\done
\endcorollary

\subsection{Real conics\pdfstr{}{ \rm(\cf.~\cite{DIS})}}\label{s.real}
Given a real structure $\rs\:X\to X$ on a $K3$-surface~$X$,
we denote by $\L_{\pm\rs}\subset\L=H_2(X;\Z)$ the
$(\pm1)$-eigenlattices of $\rs_*$.
We say that a
%real surface~$X$ or, more precisely,
real structure
$\rs\:X\to X$ is \emph{uniform} if
$\rs_*$ restricts to
$-\id\:\!\NS(X)\to\NS(X)$.
Clearly, on
%On
a uniformly real $K3$-surface~$X$, each
$(-2)$-curve is real; the converse also holds provided that $\NS(X)$ is
rationally generated by $(-2)$-curves.

Let $\rs\:X\to X$ be a real structure on a quartic $X\subset\Cp3$.
%, and let
%$\NS_{-\rs}(X)\ni h$ be the $\rs_*$-skew-invariant sublattice.
A conic $c\in\NS(X)$ is real if and only if $\rs_*(c)=-c$; hence, $\rs_*$
restricts to $-\id$ on
the primitive hull~$N$
of the
sublattice generated by real conics. Clearly, $N$ is geometric, as it is
primitive in $\NS(X)$, and
by an obvious analogue of \autoref{prop.existence}, there exists a
quartic~$X'$, arbitrarily close to~$X$, and a
uniform real structure $\rs'\:X'\to X'$ such that
\roster*
\item
the conics on~$X'$ are in a bijection with the real conics on~$X$,
\item
$\NS(X')=N$ is rationally generated by conics.
\endroster

\remark\label{rem.uniform}
Unlike~\cite{DIS}, where only lines are considered,
in general
we cannot assert that $\Fn X'$ is a subgraph of $\Fn X$: some
irreducible conics on~$X'$ may become reducible on~$X$. Nor can we assert
that $X'$ is a Barth--Bauer quartic,
even if $X$ is,
as we do not assume that $X$ has a $16$-tuple~$\dm$ of real Kummer conics.
\endremark

\lemma[\cf. {\cite[Lemma 3.8]{DIS}}]\label{lem.real}
If $\rs\:X\to X$ is a uniform real structure on a Barth--Bauer quartic
\rom(more generally, Kummer surface\rom)~$X$, then $\L_{+\rs}\cong\bU(2)$
is an orthogonal direct summand in $T(X)$.
Conversely, if
%$T(X)$ contains $\bU(2)$ as an orthogonal direct summand,
$T(X)\cong\bU(2)\oplus T'$,
then $X$ is equiconically deformation equivalent to a uniformly real quartic.
\endlemma

\proof
We have $\L_{+\rs}\subset T(X)$ and, in view of~\eqref{eq.Sperp},
$\L_{+\rs}(\frac12)\subset T(X)(\frac12)$ are also even lattices.
On the other hand,
$\discr\L_{+\rs}$ is an elementary abelian $2$-group
(see~\cite{Nikulin:forms}); hence, $\L_{+\rs}(\frac12)$
is unimodular and, as such, is
an orthogonal direct summand in any overlattice. Besides,
$\Gs_+(\L_{+\rs})=1$ (see~\cite{Nikulin:forms}) and
$\rank\L_{+\rs}\le\rank\Sperp=6$;
hence, using the classification of unimodular even lattices of small rank,
$\L_{+\rs}(\frac12)\cong\bU$.

For the converse statement, letting $N:=\NS(X)$, we mimic the arguments of
\cite{Dolgachev:polarized,Nikulin:forms}:
the period~$\Go_X$ in the
period domain of marked $N$-polarized $K3$-surfaces is moved to
$\Go_++i\Go_-$, where $\Go_+\in\bU(2)\otimes\R$ and $\Go_-\in T'\otimes\R$
are sufficiently generic vectors, $\Go_+^2=\Go_-^2>0$, and
$\R\Go_+\oplus\R\Go_-$ is positively oriented. The new quartic has the same
N\'{e}ron--Severi lattice~$N$ and, by~\cite{Nikulin:forms}
and the global Torelli theorem, it has a real
structure~$\rs$ such that $\L_{+\rs}$ is the summand~$\bU(2)$ in the
statement.
\endproof

%The following statement found in~\cite{DIS} extends literally to the case of
%conics.

\corollary\label{cor.real.singular}
If $X$ is
a Barth--Bauer quartic, the real conics on~$X$ \rom(with respect to any real
structure\rom)
span a sublattice of $\NS(X)$ or rank at most~$19$.
\done
\endcorollary

\subsection{Proof of \autoref{th.real}}\label{proof.real}
Due to \autoref{cor.real.singular}, a rigid Barth--Bauer quartic cannot have
all its conics real. However,
in view of \autoref{rem.uniform}, we cannot directly refer to
\autoref{tab.2}, and we have to compute the set of real conics
for each quartic~$X$ in \autoref{tab.3} with
$\ls|\Fn_2X|\ge560$ and each real structure $\rs\:X\to X$.
The induced action $g_N:=\rs_*$ on $N:=\NS(X)$ is found
\via~\eqref{eq.aut.group}, with an extra requirement that both
$g_N\in\OG_h(N)$ and
$g_T\in\OG(T)\sminus\OGplus(T)$,
\cf. \autoref{prop.absolute}\iref{abs.real}, must be involutive.
The computation results in a single pair $(X,\rs)$, \viz. the one marked
with a \maxreal2 in \autoref{tab.3}, and the
corresponding
graph of real conics is isomorphic to $\Fn Y$, $Y\in\famY$, where $\famY$ is
the $1$-parameter family
marked with a \maxreal1 in \autoref{tab.2}.

To complete the construction of the $1$-parameter family as in the theorem,
we
observe that $T(Y)\cong\bU(2)\oplus[76]$ for a generic $Y\in\famY$ and refer
to \autoref{lem.real}.

Finally, we check that, for each conic $c\in\Fn Y$, there is another
conic~$c'$ such that $c\cdot c'=1$, implying that $c$ has a real point.
\qed

\remark
Using the fact that, for $N:=\NS(Y)$,
the natural homomorphism $\OG_h(N)\onto\Aut(\discr N)$ is
surjective and analyzing the walls in the period domain, it is not difficult
to show that the moduli space of real quartics as in \autoref{th.real} is
connected. Technical details and analysis of the chirality are left to the
reader.
\endremark

\subsection{Proof of \autoref{prop.NS} and \autoref{cor.structure}}\label{proof.NS}
\autoref{prop.NS} is merely a combined statement of \autoref{prop.1} and
Corollaries~\ref{cor.Sh.conics},
\ref{cor.2-generated},
and~\ref{cor.3-generated}.

For \autoref{cor.structure}, we observe that, due to \autoref{prop.NS}, for
any subgraph $\dm\subset\graph$ as in the statement, the primitive hull
$(\Q\dm+\Q h)\cap\Fano(\graph)$ admits a primitive isometry to~$\L$. Hence,
it is isomorphic to~$\lattice$ as in~\eqref{eq.Sh}, which determines a unique
$4$-Kummer structure on~$\dm$ (\cf. also \autoref{lem.bS}).
\qed

{
\let\.\DOTaccent
\def\cprime{$'$}
\bibliographystyle{amsplain}
\bibliography{degt}
}

\end{document}